\documentclass[12pt,a4paper,reqno,oneside]{amsart}
\usepackage{mathabx}
\usepackage{times}

\DeclareMathAlphabet{\mathsl}{OT1}{cmss}{m}{sl}
\SetMathAlphabet{\mathsl}{bold}{OT1}{cmss}{bx}{sl}

\usepackage
{geometry}

\usepackage[dvipsnames]{xcolor}
\usepackage[english]{babel} 

\definecolor{bred}{rgb}{0.8,0,0}

\usepackage[sort,nocompress]{cite}
\usepackage{amsthm,amsmath,amssymb,bbm,geometry,epsfig,listings,paralist,enumerate,comment,nicefrac,verbatim,multirow,xcolor,wasysym,tikz}
\usepackage[pagewise]{lineno}[4.41]
\usepackage{mathrsfs}  
\usepackage{marginnote}\reversemarginpar

\usepackage[colorlinks,linkcolor=blue,citecolor=blue]{hyperref}

\usepackage[notref, notcite,final]{showkeys}
\usepackage[capitalise,compress]{cleveref} 

\crefname{equation}{}{}
\newtheorem{lemma}{Lemma}[section]

\newtheorem{proposition}[lemma]{Proposition}
\newtheorem{theorem}[lemma]{Theorem}
\newtheorem{example}[lemma]{Example}

\newtheorem{setting}[lemma]{Setting}
\crefname{subsection}{Subsection}{Subsections}
\crefname{enumi}{item}{items}
\creflabelformat{enumi}{(#2\textup{#1}#3)}

\newcommand{\1}{\ensuremath{\mathbbm{1}}}
\providecommand{\N}{{\ensuremath{\mathbbm{N}}}}
\providecommand{\Z}{{\ensuremath{\mathbbm{Z}}}}
\providecommand{\R}{{\ensuremath{\mathbbm{R}}}}

\renewcommand{\P}{{\ensuremath{\mathbbm{P}}}}
\providecommand{\E}{{\ensuremath{\mathbbm{E}}}}

\providecommand{\bfD}{{\ensuremath{\mathbf{D}}}}
\providecommand{\calR}{{\ensuremath{\mathcal{R}}}}
\providecommand{\calD}{{\ensuremath{\mathcal{D}}}}
\providecommand{\bfN}{{\ensuremath{\mathbf{N}}}}

\newcommand{\calP}{\mathcal{P}}
\newcommand{\xeqref}[1]{}
\newcommand{\threenorm}[1]{{\left\vert\kern-0.25ex\left\vert\kern-0.25ex\left\vert #1 
    \right\vert\kern-0.25ex\right\vert\kern-0.25ex\right\vert}}

\newcommand{\supnorm}[1]{{\left\vert\kern-0.25ex\left\vert\kern-0.25ex\left\vert #1 
    \right\vert\kern-0.25ex\right\vert\kern-0.25ex\right\vert}}

\newcommand{\exponentV}{{p_{\mathrm{v}}}}
\newcommand{\exponentZ}{{p_{\mathrm{z}}}}
\newcommand{\exponentX}{{p_{\mathrm{x}}}}

\newcommand{\pr}{\mathrm{pr}}

\title[]{\vspace*{-1.5cm}Rectified deep neural networks
overcome the curse of dimensionality\\
 in the numerical approximation of \\
 gradient-dependent semilinear heat equations}

\author[A. Neufeld]{Ariel Neufeld$^{1}$}
\address{$^1$  Division of Mathematical Sciences, School of Physical and Mathematical Sciences, Nanyang Technological University, Singapore}
\email{ariel.neufeld@ntu.edu.sg}

\author[T.A. Nguyen]{Tuan Anh Nguyen$^{2}$}
\address{$^2$  Faculty of Mathematics, Bielefeld University, Germany}
\email{tnguyen@math.uni-bielefeld.de}



\thanks{
Financial support by the Nanyang Assistant Professorship Grant (NAP Grant) \textit{Machine Learning based Algorithms in
Finance and Insurance} is gratefully acknowledged.}
\allowdisplaybreaks

\begin{document}
\maketitle
\begin{abstract}
\vspace*{-0.7cm}Numerical experiments indicate that deep learning algorithms overcome the curse of dimensionality when approximating solutions of semilinear 
PDEs. For certain linear PDEs and semilinear PDEs with \textit{gradient-independent} nonlinearities this has also been proved mathematically, i.e., it has been shown that
the number of parameters of the approximating
DNN increases at most polynomially in both the PDE dimension $d\in \N$ and the reciprocal of
the prescribed accuracy $\epsilon\in (0,1)$. \\
The main contribution of this paper is to rigorously prove  for the first time that deep neural networks can also overcome the curse dimensionality in the approximation of a certain class of nonlinear PDEs with \emph{gradient-dependent} nonlinearities.

\end{abstract}
\section{Introduction}

Deep learning-based approximation algorithms
for certain nonlinear parabolic partial differential equations (PDEs) have been first proposed in \cite{EHanJentzen2017,han2018solving}.
Due to their success in approximately solving nonlinear PDE in higher dimensions where classical methods have failed, there is now a variety of deep learning algorithms for  different kinds of PDEs in the scientific literature, 
inspired by the seminal work by \cite{han2018solving}.
 We refer to \cite{beck2019machine,FTT2019,GPW2022, han2020convergence, Hen2017,hure2020deep,jacquier2023deep, Rai2023} for deep BSDE methods,
to \cite{beck2020deep,beck2021deep} for deep splitting methods,
to \cite{al2022extensions, sirignano2018dgm} for deep Galerkin methods, 
to \cite{nguwi2022deep,nguwi2022numerical}  for deep branching methods,
to \cite{lu2021deepxde, raissi2019physics, zhang2020learning} based on PINNs,
to \cite{EYu2018} for the deep Ritz method,
to \cite{han2019solving} for the deep wave function method,
to \cite{BBG+2021,berner2020numerically,KLY2021,Mis2019, NM2019,nguwi2023deep} whose method is generally based on reformulating the solution of the PDE as  suitable stochastic optimization problem which can be approximately solved using some stochastic gradient descent type of algorithm,
to \cite{gonon2023random,jacquier2023random,NS2023,NSW2024} for using random neural networks,
to \cite{ito2021neural}  for solving HJB and HJBI equations related to control problems,
 and to \cite{castro2022deep,frey2022deep,gnoatto2022deep} for solving partial integro-differential equations (PIDEs).
We also refer to \cite{beck2020overview,weinan2021algorithms} for papers providing an overview on several existing deep learning  methods to solve PDEs.

Numerical experiments indicate that 
deep learning methods work exceptionally well when approximating solutions of high-dimensional PDEs and that they do not suffer from the curse of dimensionality.
However, there exist only few theoretical results proving that deep learning based approximations of solutions of PDEs do not suffer from the curse of dimensionality.

More precisely, 
\cite{BGJ2020}
 shows that empirical risk minimization  over DNN hypothesis classes overcomes the curse
of dimensionality for the numerical solution of linear Kolmogorov equations with affine coefficients.
Next,
\cite{EGJS2022} considers the pricing problem  of a European
best-of-call option on a basket of $d$ assets within the Black–Scholes model  and  proves that the solution to the $d$-variate option pricing
problem can be approximated up to an $\epsilon$-error by a deep ReLU network with depth
$\mathcal{O}(\ln(d)\ln(\epsilon^{-1}) + (\ln(d))^2)$
and $\mathcal{O}(d^{2+\frac{1}{n}}\epsilon^{-\frac{1}{n}})$
nonzero weights, where $n\in \N$ is arbitrary (with the constant implied in $\mathcal{O}(\cdot)$ depending on $n$).
Furthermore, \cite{gonon2023random}
investigates the use of random  neural networks for learning Kolmogorov PIDEs associated to Black-Scholes and more general exponential Lévy models. Here, random  neural networks are single-hidden-layer feedforward neural networks in which the input weights are randomly generated and only the output weights are trained.
Moreover, \cite{gonon2021deep} studies the expression rates of DNNs  for option
prices written on baskets of $d$ risky assets whose log-returns are modelled by a multivariate L\'evy process with general correlation structure of jumps. 
Note that the PIDEs studied by \cite{gonon2021deep} are also Black-Scholes-type PIDEs (see \cite[Display~(2.3)]{gonon2021deep}).
Next, 
\cite{gonon2023deep} proves that DNNs with ReLU activation function are able to
express viscosity solutions of Kolmogorov linear PIDEs on state spaces of
possibly high dimension $d$. Furthermore,
\cite{GHJVW2023} proves that DNNs overcome the curse of dimensionality when approximating the solutions to Black-Scholes PDEs and \cite{JSW2021} proves
that DNNs overcome the curse
of dimensionality in the numerical approximation of linear Kolmogorov PDEs with constant
diffusion and nonlinear drift coefficients. 
In addition, \cite{NS2023} proves that the  solution of the linear heat equation can be approximated by a random neural network whose amount of neurons only grow polynomially in the space dimension of the PDE and the reciprocal of the accuracy, hence overcoming the curse of dimensionality.
Moreover,
\cite{HJKN2020a}
proves that DNNs overcome the curse of dimensionality 
in the numerical approximation of semilinear heat equations and \cite{AJK+2023} extends \cite{HJKN2020a} to estimates with respect to $L^p$-norms, $p\geq 2$.
Next, \cite{CHW2022} extends \cite{HJKN2020a} to semilinear PDEs with general drift and diffusion coefficients
and \cite{neufeldnguyenwu2023deep} extends \cite{HJKN2020a} to semilinear PIDEs.

However, we highlight that all the above results \cite{HJKN2020a,BGJ2020,AJK+2023,CHW2022, neufeldnguyenwu2023deep,EGJS2022,gonon2023random,gonon2021deep,gonon2023deep,GHJVW2023,JSW2021,NS2023}
 only deal with either \emph{linear} PDEs or \emph{gradient-independent} semilinear PDEs, i.e., the case when the nonlinear part of the corresponding
semilinear PDE does not depend on the gradient of the solution. 

The contribution of our article is to prove for the first time that deep neural networks can also overcome the curse dimensionality in the approximation of nonlinear PDEs with \emph{gradient-dependent} nonlinearities. Our main result, \cref{d23}, 
proves that for semilinear heat equations
with \emph{gradient-dependent} nonlinear part  the number of parameters of the approximating
DNN increases at most polynomially in both the PDE dimension $d\in \N$ and the reciprocal of
the prescribed accuracy $\epsilon\in (0,1)$, i.e., DNNs do overcome the curse of dimensionality when approximating such PDEs.

\subsection{Notations}\label{n01}Throughout our paper we use the following notations.
Let $\left \lVert\cdot \right\rVert, \supnorm{\cdot} \colon (\cup_{d\in \N} \R^d) \to [0,\infty)$, 
$\dim \colon (\cup_{d\in \N}\R^d) \to \N$
satisfy for all $d\in \N$, $x=(x_1,\ldots,x_d)\in \R^d$ that $\|x\|=\sqrt{\sum_{i=1}^d(x_i)^2}$, $\supnorm{x}=\max_{i\in [1,d]\cap \N}|x_i|$, and
$\dim (x)=d$. Moreover, given a probability space
$(\Omega,\mathcal{F},\P)$,  a random variable $\mathfrak{X}\colon \Omega\to\R$, and
$s\in [1,\infty)$  
let 
$\lVert \mathfrak{X}\rVert_s\in [0,\infty]$ satisfy that 
$ \lVert \mathfrak{X}\rVert_s= (\E[\lvert \mathfrak{X}\rvert^s])^\frac{1}{s}$.

\subsection{A mathematical framework for DNNs}

In order to formulate our main result, \cref{d23}, we need to   introduce a mathematical framework for DNNs, see \cref{m07} below.

\begin{setting}[A mathematical framework for DNNs]\label{m07}
Let
 $\mathbf{A}_{d}\colon \R^d\to\R^d $, $d\in \N$, satisfy for all $d\in\N$, $x=(x_1,\ldots,x_d)\in \R^d$ that 
\begin{align}
\mathbf{A}_{d}(x)= \left(\max\{x_1,0\},\max\{x_2,0\},\ldots,\max\{x_d,0\}\right).
\end{align}
Let $\mathbf{D}=\cup_{H\in \N} \N^{H+2}$.
Let
\begin{align}
\mathbf{N}= \bigcup_{H\in  \N}\bigcup_{(k_0,k_1,\ldots,k_{H+1})\in \N^{H+2}}
\left[ \prod_{n=1}^{H+1} \left(\R^{k_{n}\times k_{n-1}} \times\R^{k_{n}}\right)\right].
\end{align} Let $\mathcal{D}\colon \mathbf{N}\to\mathbf{D}$, 
$\mathcal{P}\colon \mathbf{N}\to \N$,
$
\mathcal{R}\colon \mathbf{N}\to (\cup_{k,l\in \N} C(\R^k,\R^l))$
satisfy that
for all $H\in \N$, $k_0,k_1,\ldots,k_H,k_{H+1}\in \N$,
$
\Phi = ((W_1,B_1),\ldots,(W_{H+1},B_{H+1}))\in \prod_{n=1}^{H+1} \left(\R^{k_n\times k_{n-1}} \times\R^{k_n}\right), 
$
$x_0 \in \R^{k_0},\ldots,x_{H}\in \R^{k_{H}}$ with the property that
$\forall\, n\in \N\cap [1,H]\colon x_n = \mathbf{A}_{k_n}(W_n x_{n-1}+B_n )
$ we have that
\begin{align}
\mathcal{P}(\Phi)=\sum_{n=1}^{H+1}k_n(k_{n-1}+1),
\quad 
\mathcal{D}(\Phi)= (k_0,k_1,\ldots,k_{H}, k_{H+1}),
\end{align}
$
\mathcal{R}(\Phi )\in C(\R^{k_0},\R ^ {k_{H+1}}),
$
and
\begin{align}
 (\mathcal{R}(\Phi)) (x_0) = W_{H+1}x_{H}+B_{H+1}.
\end{align}
\end{setting}
Let us comment on the mathematical objects  in  \cref{m07}. For all $ d\in \N $, $ \mathbf{A}_d\colon\R^d\to\R^d$  refers to the componentwise  rectified linear unit (ReLU) activation function. 
 By $ \mathbf{N} $ we denote the set of all
(parameters characterizing) artificial feed-forward DNNs, by $ \calR $   we denote the operator that maps each DNN to its corresponding function, by $ \calP $ we denote the function that maps a DNN to its number of parameters, and by $ \calD $ we denote the function that maps a DNN to the vector of its layer dimensions. 
\subsection{Main result}
We are now in a position to state the main result, \cref{d23} below.
\begin{theorem}\label{d23}
Assume \cref{m07}. Let $T\in (0,\infty)$, $\beta,c\in[2,\infty) $, $q\in [1,2)$. For every $d\in \N$ let $L_i^d\in \R$, $i\in [0,d]\cap\Z$, satisfy that
$\sum_{i=0}^{d}L_i^d\leq c$.
For every $d\in \N$
let
 $\Lambda^d=(\Lambda^d_{\nu})_{\nu\in [0,d]\cap\Z}\colon [0,T]\to \R^{1+d}$ satisfy for all $t\in [0,T]$ that $\Lambda^d(t)=(1,\sqrt{t},\ldots,\sqrt{t})$. 
For every $d\in \N$ let $\pr^d=(\pr^d_\nu)_{\nu\in [0,d]\cap\Z}\colon \R^{d+1}\to\R$ satisfy
for all 
$w=(w_\nu)_{\nu\in [0,d]\cap\Z}$,
$i\in [0,d]\cap\Z$ that
$\pr^d_i(w)=w_i$. For every $d\in \N$, $\varepsilon\in (0,1)$ let 
$g^d,g^d_\varepsilon\in C(\R^d,\R)$, 
$f^d,f^d_\varepsilon\in C(\R^{d+1},\R)$,
$\Phi_{g^d_\varepsilon}, \Phi_{f^d_\varepsilon}\in \bfN$ satisfy that $\calR (\Phi_{g^d_\varepsilon})=g^d_\varepsilon$ and
$\calR (\Phi_{f^d_\varepsilon})=f^d_\varepsilon$.
Assume for all $d\in \N$, $\varepsilon\in (0,1)$, $x,y\in \R^d$, $w,w_1,w_2\in \R^{d+1}$ that
\begin{align}\label{b04bb}
\lvert g_\varepsilon^d(x)\rvert+\lvert Tf^d_\varepsilon(0)\rvert\leq c(d^c+\lVert x\rVert^2)^\frac{1}{2},
\end{align}
\begin{align}
\lvert
f_\varepsilon^d(w_1)-f_\varepsilon^d(w_2)\rvert\leq \sum_{\nu=0}^{d}\left[
L_\nu^d\Lambda_\nu^d(T)
\lvert\pr_\nu^d(w_1-w_2) \rvert\right]
,\label{c01bb}
\end{align}
\begin{align}
\lvert
g^d_\varepsilon(x)-g^d_\varepsilon(y)\rvert
\leq cd^c
\frac{\lVert x-y\rVert}{\sqrt{T}},\label{b19b}
\end{align}
\begin{align}
\left\lvert
{f}_\varepsilon^d(w)-f^d(w)\right\rvert\leq \frac{\varepsilon cd^c}{T}\left(1+\sum_{\nu=0}^{d}(\Lambda_\nu^d(T))^q\left\lvert\pr_\nu^d(w)\right\rvert^q\right),\label{k03b}
\end{align}
\begin{align}\label{k04b}
\left\lvert
g^d(x)-{g}^d_\varepsilon(x)\right\rvert\leq \varepsilon cd^c(d^c+\lVert x\rVert^2)^\beta,
\end{align}
\begin{align}\label{d04}
\max \!\left\{
\calP(\Phi_{g^d_\varepsilon}), 
\calP(\Phi_{f^d_\varepsilon})\right\}\leq cd^c \varepsilon^{-c}.
\end{align}
Then the following items hold.
\begin{enumerate}[(i)]
\item \label{k36} For all $d\in \N$ there exists a unique continuous function $u^d\colon [0,T)\times \R^d\to \R^{d+1}$ such that $v^d:=\pr_0^d(u^d)$ is the unique viscosity solution to the following 
semilinear PDE of parabolic type:
\begin{align}
&
\frac{\partial v^d}{\partial t}(t,x)
+\frac{1}{2}(\Delta v^d)(t,x)+f^d(v^d(t,x), (\nabla_xv^d)(t,x))
=0
\quad 
\forall\, t\in (0,T), x\in \R^d,\label{c36}
\\
&
v^d(T,x)=g^d(x) \quad \forall\, x\in \R^d,\label{c38}
\end{align}
$\nabla_xv^d=(\pr^d_1(u^d), \pr^d_2(u^d),\ldots, \pr^d_d(u^d))$,
and
\begin{align}
\max_{\nu\in [0,d]\cap\Z}\sup_{\tau\in [0,T), \xi\in \R^d}
\left[\Lambda^d_\nu(T-\tau)\frac{\lvert\pr^d_\nu(u^d(\tau,\xi))\rvert}{(1+\lVert \xi\rVert^2)^\frac{1}{2}}\right]<\infty.
\end{align}

\item There exists $\eta\in (0,\infty)$ such that for all $d\in \N$, $\epsilon\in (0,1)$ there exists $\Psi_{d,\epsilon}\in \bfN$ satisfying $\calR(\Psi_{d,\epsilon})\in C(\R^d,\R^{d+1})$, $ \calP(\Psi_{d,\epsilon})\leq \eta d^\eta\epsilon^{-\eta}$, and
\begin{align}
\left(\int_{[0,1]^d}\sum_{\nu=0}^{d}
\left\lvert\Lambda^d_\nu(T)
\pr^d_\nu\!\left(
(\calR(\Psi_{d,\epsilon}))(x)
-u^d(0,x)\right)
\right\rvert^2dx
\right)^\frac{1}{2}
\leq \epsilon.
\end{align}
\end{enumerate}
\end{theorem}
Let us make some comments on the mathematical objects in \cref{d23} above.
First, \eqref{b04bb}--\eqref{b19b} are growth and Lipschitz conditions.
Condition \eqref{k03b}--\eqref{k04b} ensure that the input functions 
$f^d$, $g^d$
can be approximated by the functions
$f^d_\varepsilon$, $g^d_\varepsilon$. The bound $cd^c\varepsilon^{-c}$ in 
condition
\eqref{d04}, which is a polynomial of $d$ and $\varepsilon^{-1}$, 
ensures that the
functions $f^d_\varepsilon$, $g^d_\varepsilon$ can be represented by DNNs without curse of dimensionality.
Note that in the \emph{gradient-independent} case \cite{HJKN2020a} we do not need to assume that $f^d$ can be approximated by DNNs without curse of dimensionality because in the one-dimensional case (cf. \cite[Corollary~3.13]{HJKN2020a}) such approximations always exist.

Under these assumptions
\cref{d23} states that, roughly speaking, if DNNs can approximate the initial condition, the linear part, and the
nonlinear part in \eqref{c36}--\eqref{c38}  without  curse of dimensionality, then they can also
approximate its solution without  curse of dimensionality.

\subsection{Outline of the proof and organization of the paper}
The proof of \cref{d23} above relies on full history
recursive multilevel Picard (MLP) approximations which have been proved to overcome the
curse of dimensionality when approximating solutions of semilinear heat equations in the gradient-dependent case, see  
\cite{HJK2022,NNW2023}. 
It is the main reason why we need to assume \eqref{c01bb}:
it ensures that MLP approximations overcome the curse of dimensionality when approximating the solution to \eqref{c36}--\eqref{c38} (cf. \cite[Theorem~1.1]{HJK2022} especially condition (1) together with the condition 
$\sum_{i=1}^d L_{d,i}\leq \lambda$ in there).
Our proof essentially consists of two main steps: \cref{c02} shows that realizations of
certain suitable MLP approximations can be represented by DNNs and \cref{b37} is a perturbation result for solutions of nonlinear PDEs of the form~\eqref{c36}.

Let us sketch the proof of \cref{d23}. First, the solution to the PDE \eqref{c36} and its gradient can be represented as solution $u^d$ to the stochastic fixed point equation (SFPE) \eqref{d11}. Next, we approximate $u^d$ by $u^d_\varepsilon$ which is the solution to SFPE \eqref{d11b} with coefficients $f^d_\varepsilon$ and $g^d_\varepsilon$. Furthermore, $u^d_\varepsilon$ can be approximated by the MLP approximations $U^{d,\theta}_{n,m,\varepsilon}$ in \eqref{d02b}. Next, we decompose the error $U^{d,\theta}_{n,m,\varepsilon}-u^d $ as follows:
$U^{d,\theta}_{n,m,\varepsilon}-u^d =(U^{d,\theta}_{n,m,\varepsilon}-u^d_\varepsilon)+(u^d_\varepsilon -u^d)$.
We estimate
$U^{d,\theta}_{n,m,\varepsilon}-u^d_\varepsilon$
using
\cite[Lemma 4.3]{NNW2023} (see  \eqref{d12}) and estimate
$u^d_\varepsilon -u^d$ using the perturbation result, \cref{b37} (see \eqref{d13}). We use \cref{c02} to show that 
$U^{d,\theta}_{n,m,\varepsilon}$ can be represented by a DNN as well as to estimate the number of layers and the maximum norm of the vector of layer dimensions needed. 
In summary, the strategy of the proof of \cref{d23} is similar to the proofs for analog results for the gradient-independent case, namely to achieve a perturbation result for the solution of the PDE under consideration together with proving that suitable MLP approximations can be represented by DNNs.  However, we highlight that the derivation of such results is much more mathematically involved for the gradient-dependent case, as the corresponding solution of the SFPE is  $(d+1)$-dimensional, instead of only one-dimensional in the gradient-independent case, due to the presence of the gradient of the solution of the PDE in the solution of the SFPE. More technically speaking, for example to prove the perturbation result \cref{b37}, we need to apply a new type of Gr\"onwall-type lemma we derived in \cite[Corollary~2.4]{NNW2023} that can deal with \eqref{k05b} and \eqref{k05c}.

The remaining part of our paper is organized as follows.
In \cref{s04b}
we prove that multilevel Picard approximations can be represented by DNNs and we provide
bounds for the number of parameters of the representing DNN.
In \cref{b37b} we provide 
the perturbation lemma, \cref{b37},
 demonstrating that the solution to the SFPE corresponding to the PDE \eqref{c36}--\eqref{c38} and its gradient are stable against
perturbations in the nonlinearity $f$ and the terminal condition $g$. The main theorem is proven in  \cref{s05}. In \cref{s02}
we give an example of
$(f^d_\varepsilon)_{d\in\N,\varepsilon\in (0,1)}$,
$(f^d)_{d\in\N}$,
$(g^d_\varepsilon)_{d\in\N,\varepsilon\in (0,1)}$,
$(g^d)_{d\in\N}$,
$(\Phi_{f^d_\varepsilon})_{d\in\N,\varepsilon\in (0,1)},(\Phi_{g^d_\varepsilon})_{d\in\N,\varepsilon\in (0,1)}$
that satisfy
the assumption of \cref{d23}.

\section{Deep neural networks}\label{s04b}
\subsection{Properties of operations associated to DNNs}\label{s04}
In \cref{m07b}
below we introduce operations which
are important for constructing the random DNN that represents the MLP approximations in the proof of
\cref{c02}.
\begin{setting}\label{m07b}
Assume \cref{m07}. Let $\odot \colon \mathbf{D}\times \mathbf{D} \to\mathbf{D} $ satisfy 
for all $H_1,H_2\in \N$, $ \alpha=(\alpha_0,\alpha_1,\ldots,\alpha_{H_1},\alpha_{H_1+1})\in\N^{H_1+2}$, $\beta=(\beta_0,\beta_1,\ldots,\beta_{H_2},\beta_{H_2+1})\in\N^{H_2+2}$
that
$
\alpha\odot \beta= (\beta_{0},\beta_{1},\ldots,\beta_{H_2},\beta_{H_2+1}+\alpha_{0},\alpha_{1},\alpha_{2},\ldots,\alpha_{H_1+1})\in \N^{H_1+H_2+3}.
$ 
Let $\boxplus \colon \mathbf{D}\times \mathbf{D} \to\mathbf{D}  $ satisfy
for all $H\in \N$, 
$\alpha= (\alpha_0,\alpha_1,\ldots,\alpha_{H},\alpha_{H+1})\in \N^{H+2}$,
$\beta= (\beta_0,\beta_1,\beta_2,\ldots,\beta_{H},\beta_{H+1})\in \N^{H+2}$
that
$
\alpha \boxplus \beta =(\alpha_0,\alpha_1+\beta_1,\ldots,\alpha_{H}+\beta_{H},\beta_{H+1})\in \N^{H+2}.
$
Let $\mathfrak{n}_n
\in \mathbf{D} $, $n\in [3,\infty)\cap\Z$,   satisfy for all $n\in [3,\infty)\cap\N$ that
\begin{align}
\label{k07}
\mathfrak{n}_n= (1,\underbrace{2,\ldots,2}_{(n-2)\text{ times}},1)\in \N^n
.
\end{align} 
 For every $n\in \N$ let $R_n\colon \bfD\to \bfD$ satisfy 
for all $H\in \N$, $\alpha=(\alpha_0,\alpha_1,\ldots,\alpha_H,\alpha_{H+1})\in \N^{H+2}$
that\begin{align}\label{k08}
 R_n(\alpha)= (\alpha_0,\alpha_1,\ldots,\alpha_H,n).
\end{align}

\end{setting}

For the proof of our main result in this section, \cref{c02}, we need several auxiliary results, Lemmas \ref{k43}--\ref{b03}, which are basic facts on DNNs. The proof of Lemmas \ref{k43}--\ref{b03}
(except \cref{b15b,p01b})
 can be found in 
\cite{HJKN2020a, CHW2022} and therefore omitted.
\cref{k43,k43b} ensure that $\odot$ and $\boxplus$ are associative so that we could write $\alpha\odot\beta\odot \gamma $
and 
$\alpha\boxplus\beta\boxplus \gamma $.
\cref{b15,b15b} are later important to estimate the maximum norm of the vector of layer dimensions of DNNs.
\cref{p01} shows that the multiplication of a DNN with a scalar can be represented by a DNN with the same vector of layer dimensions.
\cref{p01b} shows that the multiplication of a DNN with a vector  can be represented by a DNN with a slightly modified vector of layer dimensions.
\cref{m11b} shows that composition of DNNs functions can also be represented by a DNN.
\cref{b01b} shows that the sum of DNNs with the same length can also be represented by DNN. \cref{b03} is important to represent sums of DNNs which do not have the same length. 
Compared with the gradient-independent case
\cref{b15b,p01b}  are new and important in the gradient-dependent case because in the MLP approximation \eqref{h14} we need to multiply with vectors instead of only with scalars.
\begin{lemma}[$\odot$ is associative--{\cite[Lemma 3.3]{HJKN2020a}}]\label{k43}Assume \cref{m07b} and let $\alpha,\beta,\gamma\in \bfD$. Then we have that
$(\alpha\odot\beta)\odot \gamma
= \alpha\odot(\beta\odot \gamma)$.
\end{lemma}

\begin{lemma}[$\boxplus$ and associativity--{\cite[Lemma 3.4]{HJKN2020a}}]\label{k43b}Assume \cref{m07b},
let $H,k,l \in \N$, and let $\alpha,\beta,\gamma\in \left( \{k\}\times \N^{H} \times \{l\}\right)$.
Then
\begin{enumerate}[(i)]
\item we have that $\alpha\boxplus\beta\in \left(\{k\}\times \N^{H} \times \{l\}\right)$,
\item we have that $\beta\boxplus \gamma\in \left(\{k\}\times \N^{H} \times \{l\}\right)$, and 
\item we have that $(\alpha\boxplus\beta)\boxplus \gamma
= \alpha\boxplus(\beta\boxplus \gamma)$.
\end{enumerate}
\end{lemma}
 \begin{lemma}[Triangle inequality--{\cite[Lemma 3.5]{HJKN2020a}}]\label{b15}
Assume \cref{m07b},
let $H,k,l \in \N$, and let $\alpha,\beta\in \{k\}\times \N^{H} \times \{l\}$.
Then we have that
$\supnorm{\alpha\boxplus\beta}\leq\supnorm{\alpha}+
\supnorm{\beta} $.
\end{lemma}

\begin{lemma}\label{b15b}
Assume \cref{m07b},
let $H,k,l,n,m \in \N$, and let $\alpha_1,\alpha_2,\ldots,\alpha_m\in \{k\}\times \N^{H} \times \{l\}$.
Then we have that
$\supnorm{\boxplus_{i=1}^m R_n(\alpha_i) }\leq
\max\!
\left\{\sum_{i=1}^{m}
\supnorm{\alpha_i}
 ,n\right\}$.
\end{lemma}
\begin{proof}[Proof of \cref{b15b}]
Throughout this proof
for every $i\in [1,m]\cap\Z$
 let $\alpha_{i,j}\in \N$, $j\in [1,H]\cap\Z$, satisfy that
$\alpha_i=(k,\alpha_{i,1},\alpha_{i,2},\ldots,\alpha_{i,H},l)$.
 Then the definition of $\boxplus$, the definition of $R_n$, and the triangle inequality show that
$\boxplus_{i=1}^m R_n(\alpha_i) =
(k,\sum_{i=1}^m \alpha_{i,1},\sum_{i=1}^m \alpha_{i,2},\ldots,\sum_{i=1}^m \alpha_{i,H},n)
$ and
\begin{align}
\supnorm{\operatorname*{\boxplus}_{i=1}^m R_n(\alpha_i) }
&=\sup\!\left \{k,\left\lvert\sum_{i=1}^m \alpha_{i,1}\right\rvert,\left\lvert\sum_{i=1}^m \alpha_{i,2}\right\rvert,\ldots,\left\lvert\sum_{i=1}^m \alpha_{i,H}\right\rvert,n\right\}\nonumber\\
&\leq \sup\!\left \{k,\sum_{i=1}^m\left\lvert \alpha_{i,1}\right\rvert,\sum_{i=1}^m\left\lvert \alpha_{i,2}\right\rvert,\ldots,\sum_{i=1}^m\left\lvert \alpha_{i,H}\right\rvert,n\right\}\nonumber\\
&\leq \max\!
\left\{\sum_{i=1}^{m}
\supnorm{\alpha_i}
 ,n\right\}.
\end{align}
This completes the proof of \cref{b15b}.
\end{proof}

\begin{lemma}[DNNs for affine transformations--{\cite[Lemma 3.7]{HJKN2020a}}]\label{p01}
Assume \cref{m07} and let $d,m\in \N$, 
$\lambda\in \R$,
$b\in\R^d$, $a\in\R^m$, $\Psi\in\mathbf{N}$ satisfy that $\mathcal{R}(\Psi)\in C(\R^d,\R^m)$. Then we have that
$
\lambda\left((\mathcal{R}(\Psi))(\cdot +b)+a\right)\in \mathcal{R}(\{\Phi\in\mathbf{N}\colon \mathcal{D}(\Phi)=\mathcal{D}(\Psi)\}).
$
\end{lemma}

\begin{lemma}[DNNs for the multiplication with a vector]\label{p01b}
Assume \cref{m07b} and let $d,m\in \N$, 
$\lambda\in \R^m=\R^{m\times 1}$,
$a\in\R$,
$b\in\R^d$,  $\Psi\in\mathbf{N}$ satisfy that $\mathcal{R}(\Psi)\in C(\R^d,\R)$. Then we have that
$\lambda
\left((\mathcal{R}(\Psi))(\cdot +b)+a\right)\in \mathcal{R}(\{\Phi\in\mathbf{N}\colon \mathcal{D}(\Phi)= R_m(\mathcal{D}(\Psi))\}).
$
\end{lemma}

\begin{proof}
[Proof of \cref{p01b}]
Throughout this proof let $H,k_0,k_1,\ldots,k_{H+1}\in\N$ satisfy that 
\begin{align}\label{a07}
H+2=\dim(\calD(\Psi))  \quad\text{and}\quad (k_0,k_1,\ldots,k_{H},k_{H+1}) = \calD(\Psi),
\end{align}
and
let $((W_1,B_1),(W_2,B_2),\ldots,(W_H,B_H),(W_{H+1},B_{H+1})) \in \prod_{n=1}^{H+1}\left(\R^{k_n\times k_{n-1}}\times \R^{k_n}\right)$ satisfy that
\begin{align}
\Big((W_1,B_1),(W_2,B_2),\ldots,(W_H,B_H),(W_{H+1},B_{H+1})\Big)=\Psi.
\end{align}
Then the fact that $\calR(\Psi)\in C(\R^d,\R)$ implies that $k_0=d$ and $k_{H+1}=1$.
Next,
let $\phi\in \prod_{n=1}^{H}\left(\R^{k_n\times k_{n-1}}\times \R^{k_n}\right)\times (\R^{m\times k_H}\times\R^m)$ satisfy that 
\begin{align}
\phi=\Big((W_1,B_1+W_1b),(W_2,B_2),\ldots,(W_H,B_H),(\lambda W_{H+1},\lambda B_{H+1}+\lambda a)\Big),
\end{align}
This, the fact that $k_0=d$, and \eqref{a07} show 
that $\phi\in \bfN$ and 
\begin{align}
\calD(\phi)=(d,k_1,k_2,\ldots,k_H,m)= R_m(\calD(\Psi)).
\label{a10}
\end{align}
Let
$x_0,y_0 \in \R^{k_0},x_1,y_1 \in \R^{k_1},\ldots,x_{H},y_H\in \R^{k_{H}}$ satisfy 
for all $n\in \N\cap [1,H]$
that
\begin{align} 
x_n = \mathbf{A}_{k_n}(W_n x_{n-1}+B_n ),\, 
y_n = \mathbf{A}_{k_n}(W_n y_{n-1}+B_n+\1_{\{1\}}(n)W_1b ),
\quad\text{and} \quad x_0=y_0+b.
\end{align}
Then
\begin{align}
y_1= \mathbf{A}_{k_1}(W_1 y_{0}+B_1+W_1b )= \mathbf{A}_{k_1}(W_1( y_{0}+b)+B_1 )
=\mathbf{A}_{k_1}(W_1x_0+B_1 )=x_1.
\end{align}
This and an induction argument prove for all $ i\in [2,H]\cap\N$ that
\begin{align}\begin{split}
y_i=\mathbf{A}_{k_i}(W_i y_{i-1}+B_i )= \mathbf{A}_{k_i}(W_i x_{i-1}+B_i )=x_i.
\end{split}
\end{align}
This and the definition of $\calR$   prove that
\begin{align}\begin{split}
(\calR(\phi))(y_0)&= \lambda W_{H+1}y_H+\lambda B_{H+1}+\lambda a\\&=\lambda W_{H+1}x_H+\lambda B_{H+1}+\lambda a\\&
=\lambda (W_{H+1}x_H+ B_{H+1}+ a)
\\&
=\lambda(
(\calR(\Psi))(x_0)+a)\\&= \lambda (\calR(\Psi))(y_0+b)+a.
\end{split}
\end{align}
This and the fact that
$y_0$ was arbitrary
prove that 
\begin{align}
\calR(\phi)=\lambda ((\calR(\Psi))(\cdot+b)+a). 
\end{align}
This and \eqref{a10} complete the proof of \cref{p01b}. 
\end{proof}

\begin{lemma}[Composition of functions generated by DNNs--{\cite[Lemma 3.8]{HJKN2020a}}]\label{m11b}
Assume~\cref{m07b} and let $d_1,d_2,d_3\in\N$, $f\in C(\R^{d_2},\R^{d_3})$, $g\in C(  \R^{d_1}, \R^{d_2}) $, 
$\alpha,\beta\in \mathbf{D}$ satisfy that
$f\in \mathcal{R}(\{\Phi\in \mathbf{N}\colon \mathcal{D}(\Phi)=\alpha\})$
and
$g\in \mathcal{R}(\{\Phi\in \mathbf{N}\colon \mathcal{D}(\Phi)=\beta\})$.
Then we have
that $(f\circ g)\in \mathcal{R}(\{\Phi\in \mathbf{N}\colon \mathcal{D}(\Phi)=\alpha\odot\beta\})$.
\end{lemma}

\begin{lemma}[Sum of DNNs of the same length--{\cite[Lemma 3.9]{HJKN2020a}}]
\label{b01b}
Assume \cref{m07b} and let $M,H,p,q\in \N$,  $h_1,h_2,\ldots,h_M\in\R$,
 $k_i\in \mathbf{D} $,
$f_i\in C(\R^{p},\R^{q})$,
$i\in [1,M]\cap\N$, satisfy 
for all $i\in [1,M]\cap\N$
that $ \dim(k_i)=H+2$ and
$f_i\in 
\mathcal{R}(\{\Phi\in\mathbf{N}\colon \mathcal{D}(\Phi)=k_i\}).
$
Then
we have that 
$
\sum_{i=1}^{M}h_if_i
\in\mathcal{R}\left(\left\{ \Phi\in\mathbf{N}\colon
\mathcal{D}(\Phi)=\boxplus_{i=1}^Mk_i\right\}\right).
$
\end{lemma}

\begin{lemma}[Existence of DNNs with $H$ hidden layers for $\mathrm{Id}_{\R}$--{\cite[Lemma 3.6]{HJKN2020a}}]\label{b03}
Assume \cref{m07b} and let $H\in \N$.
Then we have that
$\mathrm{Id}_{\R}\in \mathcal{R}(\{\Phi\in\mathbf{N}\colon\mathcal{D}(\Phi)=\mathfrak{n}_{H+2} \}) $.
\end{lemma}


\subsection{DNN representation of  MLP approximations}In \cref{c02} below we prove that the MLP approximations under consideration can be represented by DNNs.
\begin{proposition}\label{c02}
Assume \cref{m07}.
Let 
$ T \in (0,\infty) $, 
$ d,m \in \N $, 
$
\Theta = \cup_{ n \in \N } \Z^n
$. Let $c\in (0,\infty)$ satisfy that 
\begin{align}
c\geq \max \{d+1,\supnorm{\calD (\Phi_f)}, \supnorm{\calD (\Phi_g)}\}.\label{h15}
\end{align}
Let $\varrho \colon[0,T]^2\to\R$
satisfy for all $t\in [0,T)$, $s\in (t,T)$ that
\begin{align}
\varrho(t,s)=\frac{1}{\mathrm{B}(\tfrac{1}{2},\tfrac{1}{2})}\frac{1}{\sqrt{(T-s)(s-t)}}.
\end{align}
Let $ ( \Omega, \mathcal{F}, \P )$
be a probability space.
Let
$
  W^{ \theta }\colon [0,T] \times \Omega \to \R^d 
$, 
$ \theta \in \Theta $,
be standard
$(\mathbb{F}_t)_{t\in[0,T]}$-Brownian motions
with continuous sample paths.
Let $\mathfrak{t}^\theta\colon \Omega\to(0,1)$, $\theta\in \Theta$, be 
i.i.d.\ random variables. Assume
for all $b\in (0,1)$
that 
\begin{align}
\P(\mathfrak{t}^0\le b)=\frac{1}{\mathrm{B}(\tfrac{1}{2},\tfrac{1}{2})}\int_0^b \frac{dr}{\sqrt{r(1-r)}}.
\end{align} 
Assume that
$(W^\theta)_{\theta \in \Theta}$ and $(\mathfrak{t}^\theta)_{\theta \in \Theta}$
are independent. Let $\mathfrak{T}^{\theta}\colon [0,T)\times \Omega\to [0,T)$, $\theta\in\Theta$,
 satisfy for all  $t\in [0,T)$
that $\mathfrak{T}^{\theta} _t = t+ (T-t)\mathfrak{t}^{\theta}$.
Let
$ f\in C([0,T]\times\R^d\times\R^{1+d},\R)$, 
$ 
  g\in C(\R^d, \R)
$. 
Let
$ 
  {U}_{ n,m}^{\theta }
  \colon[0,T)\times\R^d\times\Omega\to\R^{1+d}
$,
$n\in\Z$, $\theta\in\Theta$, satisfy
for all 
$
  n \in \N
$,
$ \theta \in \Theta $,
$ t\in [0,T)$,
$x \in \R^d$
that $
{U}_{-1,m}^{\theta}(t,x)={U}_{0,m}^{\theta}(t,x)=0$ and
\begin{equation}  \begin{split}
&{U}_{n,m}^{\theta}(t,x)
=(g(x), 0  ) + \sum_{i=1}^{m^n}
\left(1,\tfrac{W^{(\theta, 0, -i)}_{T}- W^{(\theta, 0, -i)}_{t}
  }{ T - t }
  \right)^\top
\tfrac{\big(g(x+W^{(\theta,0,-i)}_T-W^{(\theta,0,-i)}_t)-g(x)\big)}{m^n}
\\
&+\sum_{\ell=0}^{n-1}\sum_{i=1}^{m^{n-\ell}}
 \left(
  1 ,
  \tfrac{ 
  W_{\mathfrak{T}^{(\theta, \ell,i)}_t}^{(\theta,\ell,i)}- W^{(\theta,\ell, i)}_{t}
  }{ \mathfrak{T}^{(\theta, \ell,i)}_t-t}
  \right)^\top
\tfrac{  \left(f\circ {U}_{\ell,m}^{(\theta,\ell,i)}-\1_{\N}(\ell)f\circ  {U}_{\ell-1,m}^{(\theta,-\ell,i)}\right)
(\mathfrak{T}^{(\theta,\ell,i)}_t,x+W_{\mathfrak{T}^{(\theta, \ell,i)}_t}^{(\theta,\ell,i)}-W_t^{(\theta,\ell,i)})} {m^{n-\ell}\varrho(t,\mathfrak{T}^{(\theta, \ell,i)}_t)}
,
\end{split}  \label{h14}   \end{equation}
and let $\omega\in \Omega$.
Then for all $ n\in \N_0$ there exists $(\Phi^\theta_{n,t})_{\theta\in \Theta, t\in [0,T)}\subseteq \bfN$ such that for all $t_1,t_2,t\in [0,T)$, $\theta_1,\theta_2,\theta\in \Theta$ we have that
\begin{align}
\calD(\Phi^{\theta_1}_{n,t_1})=
\calD(\Phi^{\theta_2}_{n,t_2}),
\end{align} 
\begin{align}
\dim (\calD(\Phi^\theta_{n,t}))= n (\dim (\calD (\Phi_f))-1 )
+\dim (\calD(\Phi_g)),
\end{align}
\begin{align}
\supnorm{\calD(\Phi^\theta_{n,t})}\leq c(4m)^n,
\end{align}
\begin{align}
{U}_{n,m}^{\theta}(t,x,\omega)= (\calR(\Phi^\theta_{n,t}))(x).
\end{align}
\end{proposition}
\begin{proof}[Proof of \cref{c02}]We prove the lemma by induction on $n\in \N$. Since the zero function can be represented by a DNN of arbitrary length, the base case $n=0$ is clear.
For the induction step 
from $n\in \N_0$
to $n+1\in \N$ let $n\in \N_0$, assume for all $\ell\in [0,n]\cap\Z$,
$t_1,t_2,t\in [0,T)$, $\theta_1,\theta_2,\theta\in \Theta$
that
\begin{align}\label{h03}
\calD(\Phi^{\theta_1}_{\ell,t_1})=
\calD(\Phi^{\theta_2}_{\ell,t_2}),
\end{align} 
\begin{align}\label{h04}
\dim (\calD(\Phi^\theta_{\ell,t}))= \ell (\dim (\calD (\Phi_f))-1 )
+\dim (\calD(\Phi_g)),
\end{align}
\begin{align}\label{h05}
\supnorm{\calD(\Phi^\theta_{\ell,t})}\leq c(3m)^\ell,
\end{align}
\begin{align}\label{h06}
{U}_{\ell,m}^{\theta}(t,x,\omega)= (\calR(\Phi^\theta_{\ell,t}))(x),
\end{align}
and let the notations in \cref{m07b} be given.
First, \cref{p01,m11b,b03} imply for all $i\in [1,m^{n+1}]\cap\Z$, $\theta\in \Theta$, $t\in [0,T)$ that
\begin{align} \begin{split} 
&g
=\mathrm{Id}_\R \circ g
\in 
\calR\! \left(\left\{\Phi\in \bfN\colon 
\calD(\Phi)=
\mathfrak{n}_{(n+1)(\dim (\calD(\Phi_f) ) -1)+1}\odot 
\calD (\Phi_g)
\right\}\right)\end{split}
\end{align}
and\footnote{For every $d_1,d_2\in \N$, $z\in \R^{d_1}$, $h\in \R^{d_1}\to\R^{d_2}$ we denote by $h(\cdot +z)$ the function
$\R^{d_1}\ni x\mapsto h(x+z)\in \R^{d_2}$.}
\begin{align} 
&g\!\left(\cdot+W^{(\theta,0,-i)}_T(\omega)-W^{(\theta,0,-i)}_t(\omega)\right)=
(\mathrm{Id}_\R \circ g)\!\left(\cdot+W^{(\theta,0,-i)}_T(\omega)-W^{(\theta,0,-i)}_t(\omega)\right)\nonumber\\
&\in 
\calR \left(\left\{\Phi\in \bfN\colon 
\calD(\Phi)=
\mathfrak{n}_{(n+1)(\dim (\calD(\Phi_f) ) -1)+1}
\odot \calD(\Phi_g)
\right\}\right).\label{h09b}
\end{align}
This and \cref{p01b} imply for all $i\in [1,m^{n+1}]\cap\Z$, $\theta\in \Theta$, $t\in [0,T)$ that
\begin{align}
&
\left(0,\frac{W^{(\theta, 0, -i)}_{T}(\omega)- W^{(\theta, 0, -i)}_{t}(\omega)
  }{ T - t }
  \right)^\top
g(x)\nonumber\\
&\in \calR\! \left(\left\{\Phi\in \bfN\colon 
\calD(\Phi)=
R_{d+1}\!\left(
\mathfrak{n}_{(n+1)(\dim (\calD(\Phi_f) ) -1)+1}\odot 
\calD (\Phi_g)\right)
\right\}\right),\label{h09}
\end{align} and
\begin{align}
&
\left(1,\frac{W^{(\theta, 0, -i)}_{T}(\omega)- W^{(\theta, 0, -i)}_{t}(\omega)
  }{ T - t }
  \right)^\top
g\!\left(\cdot +W^{(\theta,0,-i)}_T(\omega)-W^{(\theta,0,-i)}_t(\omega)\right)\nonumber\\
&
\in 
\calR\! \left(\left\{\Phi\in \bfN\colon 
\calD(\Phi)=
R_{d+1}\!\left(
\mathfrak{n}_{(n+1)(\dim (\calD(\Phi_f) ) -1)+1}\odot 
\calD (\Phi_g)\right)
\right\}\right).\label{h10}
\end{align}
Next, the induction hypothesis (see \eqref{h03} and \eqref{h06}), the fact that $f=\calR (\Phi_f)$, and \cref{m11b} show for all $i\in [1,m]$, $\theta\in \Theta$, $t\in [0,T)$ that
\begin{align}
&
\left(f\circ {U}_{n,m}^{(\theta,n,i)}\right)\!
\left(\mathfrak{T}^{(\theta,n,i)}_t(\omega),\cdot +W_{\mathfrak{T}^{(\theta, n,i)}_t(\omega)}^{(\theta,n,i)}(\omega)-W_t^{(\theta,n,i)}(\omega),\omega\right)\nonumber\\&\in \calR
\!\left(
\left\{
\Phi\in \bfN \colon 
 \calD(\Phi_f)
\odot \calD (\Phi_{n,0}^0)
\right\}
\right).
\end{align}
This and \cref{p01b} imply for all $i\in [1,m]$, $\theta\in \Theta$, $t\in [0,T)$ that
\begin{align}
&
\left(1 ,\frac{W_{\mathfrak{T}^{(\theta, n,i)}_t(\omega)}^{(\theta,n,i)}(\omega)- W^{(\theta,n, i)}_{t}(\omega)}{\mathfrak{T}^{(\theta, n,i)}_t(\omega)-t}
  \right)^\top
\left(f\circ {U}_{n,m}^{(\theta,n,i)}\right)\!
\left(\mathfrak{T}^{(\theta,n,i)}_t(\omega),\cdot +W_{\mathfrak{T}^{(\theta, n,i)}_t(\omega)}^{(\theta,n,i)}(\omega)-W_t^{(\theta,n,i)}(\omega), \omega\right)\nonumber\\
&\in \calR
\!\left(
\left\{
\Phi\in \bfN \colon 
R_{d+1}
 \!\left(\calD(\Phi_f)
\odot \calD (\Phi_{n,0}^0)\right)
\right\}
\right).\label{h11}
\end{align}
Next, \cref{b03,m11b}, the fact that $f=\calR(\Phi_f)$, and the induction hypothesis  (see \eqref{h03} and \eqref{h06})
prove for all $\ell\in [0,n]\cap\Z$, $\theta\in \Theta$, $i\in [m^{n+1-\ell}]\cap\Z$, $t\in [0,T)$
that
\begin{align}
&
\left(f\circ {U}_{\ell,m}^{(\theta,\ell,i)}\right)\!
\left(\mathfrak{T}^{(\theta,\ell,i)}_t(\omega),\cdot +W_{\mathfrak{T}^{(\theta, \ell,i)}_t(\omega)}^{(\theta,\ell,i)}(\omega)-W_t^{(\theta,\ell,i)}(\omega),\omega\right)\nonumber\\
&
=
\left(\mathrm{Id}_{\R}\circ f\circ {U}_{\ell,m}^{(\theta,\ell,i)}\right)\!
\left(\mathfrak{T}^{(\theta,\ell,i)}_t(\omega),\cdot +W_{\mathfrak{T}^{(\theta, \ell,i)}_t(\omega)}^{(\theta,\ell,i)}(\omega)-W_t^{(\theta,\ell,i)}(\omega),\omega\right)\nonumber\\
&\in \calR
\!\left(
\left\{
\Phi\in \bfN \colon \mathfrak{n}_{(n-\ell) (\dim (\calD(\Phi_f))-1)+1 }
\odot \calD(\Phi_f)
\odot \calD (\Phi_{\ell,0}^0)
\right\}
\right).
\end{align}
This and \cref{p01b} 
demonstrate for all $\ell\in [0,n]\cap\Z$, $\theta\in \Theta$, $i\in [m^{n+1-\ell}]\cap\Z$, $t\in [0,T)$
that
\begin{align}
&\left(
  1 ,
  \frac{ 
  W_{\mathfrak{T}^{(\theta, \ell,i)}_t(\omega)}^{(\theta,\ell,i)}(\omega)- W^{(\theta,\ell, i)}_{t}(\omega)
  }{ \mathfrak{T}^{(\theta, \ell,i)}_t(\omega)-t}
  \right)^\top
  \left(f\circ {U}_{\ell,m}^{(\theta,\ell,i)}\right)\!
\left(\mathfrak{T}^{(\theta,\ell,i)}_t(\omega),\cdot +W_{\mathfrak{T}^{(\theta, \ell,i)}_t(\omega)}^{(\theta,\ell,i)}(\omega)-W_t^{(\theta,\ell,i)}(\omega),\omega\right)\nonumber\\
&\in \calR
\!\left(
\left\{
\Phi\in \bfN \colon 
R_{d+1}\!\left(
\mathfrak{n}_{(n-\ell) (\dim (\calD(\Phi_f))-1)+1 }
\odot \calD(\Phi_f)
\odot \calD (\Phi_{\ell,0}^0)
\right)
\right\}
\right).\label{h12}
\end{align}
Similarly, we have 
for all $\ell\in [0,n]\cap\Z$, $\theta\in \Theta$, $i\in [m^{n+1-\ell}]\cap\Z$, $t\in [0,T)$ that
\begin{align}
&
  \left(f\circ {U}_{\ell-1,m}^{(\theta,-\ell,i)}\right)\!
\left(\mathfrak{T}^{(\theta,\ell,i)}_t(\omega),\cdot +W_{\mathfrak{T}^{(\theta, \ell,i)}_t(\omega)}^{(\theta,\ell,i)}(\omega)-W_t^{(\theta,\ell,i)}(\omega),\omega\right)\nonumber\\
&
= \left(\mathrm{Id}_{\R}\circ f\circ {U}_{\ell-1,m}^{(\theta,-\ell,i)}\right)\!
\left(\mathfrak{T}^{(\theta,\ell,i)}_t(\omega),\cdot +W_{\mathfrak{T}^{(\theta, \ell,i)}_t(\omega)}^{(\theta,\ell,i)}(\omega)-W_t^{(\theta,\ell,i)}(\omega)\right)\nonumber\\
&\in \calR
\!\left(
\left\{
\Phi\in \bfN \colon \mathfrak{n}_{(n-\ell+1) (\dim (\calD(\Phi_f))-1)+1 }
\odot \calD(\Phi_f)
\odot \calD (\Phi_{\ell-1,0}^0)
\right\}
\right)
\end{align}
and
\begin{align}
 &\left(
  1 ,
  \frac{ 
  W_{\mathfrak{T}^{(\theta, \ell,i)}_t(\omega)}^{(\theta,\ell,i)}(\omega)- W^{(\theta,\ell, i)}_{t}(\omega)
  }{ \mathfrak{T}^{(\theta, \ell,i)}_t(\omega)-t}
  \right)^\top
  \left(f\circ {U}_{\ell-1,m}^{(\theta,-\ell,i)}\right)\!
\left(\mathfrak{T}^{(\theta,\ell,i)}_t(\omega),\cdot +W_{\mathfrak{T}^{(\theta, \ell,i)}_t(\omega)}^{(\theta,\ell,i)}(\omega)-W_t^{(\theta,\ell,i)}(\omega),\omega\right)\nonumber\\
&\in \calR
\!\left(
\left\{
\Phi\in \bfN \colon
R_{d+1}\!\left(
 \mathfrak{n}_{(n-\ell+1) (\dim (\calD(\Phi_f))-1)+1 }
\odot \calD(\Phi_f)
\odot \calD (\Phi_{\ell-1,0}^0)
\right)
\right\}
\right).\label{h13}
\end{align}
Next, \eqref{k08}, the definition of $\odot$, and \eqref{k07} imply that
\begin{align}
&\dim \!\left(
R_{d+1}\!\left(
\mathfrak{n}_{(n+1)(\dim (\calD(\Phi_f) ) -1)+1}\odot 
\calD (\Phi_g)\right)\right)\nonumber\\
&=\dim \!\left(
\mathfrak{n}_{(n+1)(\dim (\calD(\Phi_f) ) -1)+1}\odot 
\calD (\Phi_g)\right)\nonumber\\
&=(n+1)(\dim (\calD(\Phi_f) ) -1)+1+\dim(\calD (\Phi_g))-1\nonumber\\
&=(n+1)(\dim (\calD(\Phi_f) ) -1)+\dim(\calD (\Phi_g)).\label{h07}
\end{align}
Furthermore, \eqref{k08}, the definition of $\odot$, and the induction hypothesis (see \eqref{h04}) show that
\begin{align}
&
\dim \!\left(
R_{d+1}
 \!\left(\calD(\Phi_f)
\odot \calD (\Phi_{n,0}^0)\right)\right)\nonumber\\
&=
\dim \!\left(
\calD(\Phi_f)
\odot \calD (\Phi_{n,0}^0)\right)\nonumber\\
&=\dim (\calD(\Phi_f))+ \dim (\calD (\Phi_{n,0}^0))-1\nonumber\\
&=\dim (\calD(\Phi_f))+\left[n (\dim (\calD (\Phi_f))-1 )
+\dim (\calD(\Phi_g))\right]-1\nonumber\\
&=
(n+1)(\dim (\calD(\Phi_f) ) -1)+\dim(\calD (\Phi_g)).
\end{align}
Moreover, \eqref{k08}, the definition of $\odot$, \eqref{k07},  the induction hypothesis (see \eqref{h04}) imply that
\begin{align}
&
\dim \!\left(R_{d+1}\!\left(
\mathfrak{n}_{(n-\ell) (\dim (\calD(\Phi_f))-1)+1 }
\odot \calD(\Phi_f)
\odot \calD (\Phi_{\ell,0}^0)
\right)\right)\nonumber\\
&=\dim \!\left(
\mathfrak{n}_{(n-\ell) (\dim (\calD(\Phi_f))-1)+1 }
\odot \calD(\Phi_f)
\odot \calD (\Phi_{\ell,0}^0)
\right)\nonumber\\
&=(n-\ell) (\dim (\calD(\Phi_f))-1)+1 +\dim(\calD(\Phi_f))+\dim (\calD (\Phi_{\ell,0}^0))-2\nonumber\\
&=
(n-\ell) (\dim (\calD(\Phi_f))-1)+1 +\dim(\calD(\Phi_f))+\left[\ell (\dim (\calD (\Phi_f))-1 )
+\dim (\calD(\Phi_g))\right]-2\nonumber\\
&=(n+1)(\dim (\calD(\Phi_f) ) -1)+\dim(\calD (\Phi_g)).\label{h08}
\end{align}
Now, \eqref{h07}--\eqref{h08} show, roughly speaking, that the functions in \eqref{h09}, \eqref{h10}, \eqref{h11}, \eqref{h12}, and \eqref{h13} can be represented by networks with the same number of layers: $(n+1)(\dim (\calD(\Phi_f) ) -1)+\dim(\calD (\Phi_g))$. Hence, \cref{b01b}, \eqref{h09}, \eqref{h10}, \eqref{h11}, \eqref{h12}, and \eqref{h13} imply that there exists 
$(\Phi_{n+1,t}^\theta)_{\theta\in \Theta,t\in [0,T)}\subseteq\bfN$ such that for all $\theta\in \Theta$, $t\in [0,T) $, $x\in \R^d$ we have that
\begin{align}
&(\calR(\Phi_{n+1,t}^\theta))(x)\nonumber
\\
&=
\sum_{i=1}^{m^{n+1}}
\left(1,\tfrac{W^{(\theta, 0, -i)}_{T}- W^{(\theta, 0, -i)}_{t}
  }{ T - t }
  \right)^\top
\tfrac{g(x+W^{(\theta,0,-i)}_T-W^{(\theta,0,-i)}_t)}{m^{n+1}}\nonumber\\
&\quad -\sum_{i=1}^{m^{n+1}}\left(0,\tfrac{W^{(\theta, 0, -i)}_{T}- W^{(\theta, 0, -i)}_{t}
  }{ T - t }
  \right)^\top \tfrac{g(x)}{m^{n+1}}\nonumber\\
&\quad +\sum_{i=1}^{m}
\left(1 ,\tfrac{W_{\mathfrak{T}^{(\theta, n,i)}_t}^{(\theta,n,i)}- W^{(\theta,n, i)}_{t}}{\mathfrak{T}^{(\theta, n,i)}_t-t}
  \right)^\top
\tfrac{\left(f\circ {U}_{n,m}^{(\theta,n,i)}\right)\!
\left(\mathfrak{T}^{(\theta,n,i)}_t,x+W_{\mathfrak{T}^{(\theta, n,i)}_t}^{(\theta,n,i)}-W_t^{(\theta,n,i)}\right)} {m^{n+1-n}\varrho(t,\mathfrak{T}^{(\theta, n,i)}_t)}\nonumber\\
&\quad +\sum_{\ell=0}^{n-1}\sum_{i=1}^{m^{n+1-\ell}}
 \left(
  1 ,
  \tfrac{ 
  W_{\mathfrak{T}^{(\theta, \ell,i)}_t}^{(\theta,\ell,i)}- W^{(\theta,\ell, i)}_{t}
  }{ \mathfrak{T}^{(\theta, \ell,i)}_t-t}
  \right)^\top
\tfrac{  \left(f\circ {U}_{\ell,m}^{(\theta,\ell,i)}\right)\!
\left(\mathfrak{T}^{(\theta,\ell,i)}_t,x+W_{\mathfrak{T}^{(\theta, \ell,i)}_t}^{(\theta,\ell,i)}-W_t^{(\theta,\ell,i)}\right)} {m^{n+1-\ell}\varrho(t,\mathfrak{T}^{(\theta, \ell,i)}_t)}\nonumber
\\
&\quad -\sum_{\ell=1}^{n}\sum_{i=1}^{m^{n+1-\ell}}
 \left(
  1 ,
  \tfrac{ 
  W_{\mathfrak{T}^{(\theta, \ell,i)}_t}^{(\theta,\ell,i)}- W^{(\theta,\ell, i)}_{t}
  }{ \mathfrak{T}^{(\theta, \ell,i)}_t-t}
  \right)^\top
\tfrac{  \left(f\circ {U}_{\ell-1,m}^{(\theta,-\ell,i)}\right)\!
\left(\mathfrak{T}^{(\theta,\ell,i)}_t,x+W_{\mathfrak{T}^{(\theta, \ell,i)}_t}^{(\theta,\ell,i)}-W_t^{(\theta,\ell,i)}\right)} {m^{n+1-\ell}\varrho(t,\mathfrak{T}^{(\theta, \ell,i)}_t)}\nonumber\\
&=(g(x), 0  ) + \sum_{i=1}^{m^{n+1}}
\left(1,\tfrac{W^{(\theta, 0, -i)}_{T}- W^{(\theta, 0, -i)}_{t}
  }{ T - t }
  \right)^\top
\tfrac{g(x+W^{(\theta,0,-i)}_T-W^{(\theta,0,-i)}_t)-g(x)}{m^{n+1}}\nonumber
\\
&\quad +\sum_{\ell=0}^{n}\sum_{i=1}^{m^{n+1-\ell}}
 \left(
  1 ,
  \tfrac{ 
  W_{\mathfrak{T}^{(\theta, \ell,i)}_t}^{(\theta,\ell,i)}- W^{(\theta,\ell, i)}_{t}
  }{ \mathfrak{T}^{(\theta, \ell,i)}_t-t}
  \right)^\top
\tfrac{  \left(f\circ {U}_{\ell,m}^{(\theta,\ell,i)}-\1_{\N}(\ell)f\circ {U}_{\ell-1,m}^{(\theta,-\ell,i)}\right)\!
\left(\mathfrak{T}^{(\theta,\ell,i)}_t,x+W_{\mathfrak{T}^{(\theta, \ell,i)}_t}^{(\theta,\ell,i)}-W_t^{(\theta,\ell,i)}\right)} {m^{n+1-\ell}\varrho(t,\mathfrak{T}^{(\theta, \ell,i)}_t)}\nonumber\\
&=
{U}_{n+1,m}^{\theta}(t,x, \omega),\label{h18}
   \end{align}
\begin{align}
\dim(\Phi_{n+1,t}^\theta))=
(n+1)(\dim (\calD(\Phi_f) ) -1)+\dim(\calD (\Phi_g)),\label{h17}
\end{align}
and
\begin{align}
\calD(\Phi_{n+1,t}^\theta)&=\left(\operatorname*{\boxplus}_{i=1}^{m^{n+1}}R_{d+1}\!\left(
\mathfrak{n}_{(n+1)(\dim (\calD(\Phi_f) ) -1)+1}\odot 
\calD (\Phi_g)\right)\right)\nonumber\\
&\quad \boxplus \left(
\operatorname*{\boxplus}_{i=1}^{m^{n+1}}R_{d+1}\!\left(
\mathfrak{n}_{(n+1)(\dim (\calD(\Phi_f) ) -1)+1}\odot 
\calD (\Phi_g)\right)\right)\nonumber\\
&\quad \boxplus\left(
\operatorname*{\boxplus}_{i=1}^{m}
 R_{d+1}
 \!\left(\calD(\Phi_f)
\odot \calD (\Phi_{n,0}^0)\right)
\right)\nonumber\\
&\quad \boxplus \left(\operatorname*{\boxplus}_{\ell=0}^{n-1}
\operatorname*{\boxplus}_{i=1}^{m^{n+1-\ell}}R_{d+1}\!\left(
\mathfrak{n}_{(n-\ell) (\dim (\calD(\Phi_f))-1)+1 }
\odot \calD(\Phi_f)
\odot \calD (\Phi_{\ell,0}^0)
\right)
\right)\nonumber\\
&\quad \boxplus\left(\operatorname*{\boxplus}_{\ell=1}^{n}
\operatorname*{\boxplus}_{i=1}^{m^{n+1-\ell}}R_{d+1}\!\left(
 \mathfrak{n}_{(n-\ell+1) (\dim (\calD(\Phi_f))-1)+1 }
\odot \calD(\Phi_f)
\odot \calD (\Phi_{\ell-1,0}^0)
\right)
\right).
\end{align} 
This shows for all $t_1,t_2\in [0,T)$, $\theta_1,\theta_2\in \Theta$ that
\begin{align}
\calD(
\Phi_{n+1,t_1}^{\theta_1})=\calD(
\Phi_{n+1,t_2}^{\theta_2}).\label{h16}
\end{align}
Next, the induction hypothesis (see \eqref{h04}) and \cref{b15b}
imply for all $\theta\in \Theta$, $t\in [0,T)$ that
\begin{align}
\supnorm{\calD(\Phi_{n+1,t}^\theta)}&\leq \max \Biggl\{
\sum_{i=1}^{m^{n+1}}\supnorm{
\mathfrak{n}_{(n+1)(\dim (\calD(\Phi_f) ) -1)+1}\odot 
\calD (\Phi_g)}\nonumber\\
&\quad +\sum_{i=1}^{m^{n+1}}\supnorm{
\mathfrak{n}_{(n+1)(\dim (\calD(\Phi_f) ) -1)+1}\odot 
\calD (\Phi_g)}\nonumber\\
&\quad +\sum_{i=1}^{m}\supnorm{\calD(\Phi_f)
\odot \calD (\Phi_{n,0}^0)}\nonumber\\
&\quad +\sum_{\ell=0}^{n-1}\sum_{i=1}^{m^{n+1-\ell}}\supnorm{
\mathfrak{n}_{(n-\ell) (\dim (\calD(\Phi_f))-1)+1 }
\odot \calD(\Phi_f)
\odot \calD (\Phi_{\ell,0}^0)}\nonumber\\
&\quad +\sum_{\ell=1}^n\sum_{i=1}^{m^{n+1-\ell}}
 \mathfrak{n}_{(n-\ell+1) (\dim (\calD(\Phi_f))-1)+1 }
\odot \calD(\Phi_f)
\odot \calD (\Phi_{\ell-1,0}^0),d+1
\Biggr\}.
\end{align}
Note that the definition of $\odot$ show that for all $H_1,H_2,\alpha_0,\alpha_1,\ldots, \alpha_{H_1+1}, \beta_0,\beta_1,\ldots,\beta_{H_2+1}\in \N$, $\alpha,\beta\in \bfD$ with $\alpha=(\alpha_1,\ldots,\alpha_{H_1+1})$,
$\beta=(\beta_1,\ldots,\beta_{H_2+1})$,
$\alpha_0=\beta_{H_2+1}=1$
 we have that 
$\supnorm{\alpha\odot\beta}\leq \max \{\supnorm{\alpha},\supnorm{\beta},2\}$. This, \eqref{h16}, \eqref{k07}, \eqref{h15}, and the induction hypothesis (see \eqref{h05})
prove for all $\theta\in \Theta$, $t\in [0,T)$ that
\begin{align}
&
\supnorm{\calD(\Phi_{n+1,t}^\theta)}\nonumber
\\
&\leq 
2\left[\sum_{i=1}^{m^{n+1}}c\right]
+\left[\sum_{i=1}^{m}c(4m)^n\right]+
\left[
\sum_{\ell=0}^{n-1}\sum_{i=1}^{m^{n+1-\ell}}c(4m)^\ell\right]
+\left[\sum_{\ell=1}^n\sum_{i=1}^{m^{n+1-\ell}}c(4m)^{\ell-1}\right]\nonumber\\
&\leq 2cm^{n+1}+mc(4m)^n+
\left[
\sum_{\ell=0}^{n-1}m^{n+1-\ell}c(4m)^\ell\right]
+\left[\sum_{\ell=1}^nm^{n+1-\ell}c(4m)^{\ell-1}\right]\nonumber\\
&\leq 2cm^{n+1}+cm^{n+1}4^n+cm^{n+1}\left[\sum_{\ell=0}^{n-1}4^\ell\right]+
cm^{n+1}\left[\sum_{\ell=1}^{n}4^{\ell-1}\right]\nonumber\\
&\leq 
cm^{n+1}\left[2+2
\sum_{\ell=0}^{n}4^\ell
\right]
\leq 
cm^{n+1}\left[1+3\sum_{\ell=0}^{n}4^\ell
\right]
= cm^{n+1}\left[1+3\frac{4^{n+1}-1}{4-1}\right]\nonumber\\
&=c(4m)^{n+1}.
\end{align}
This, \eqref{h16}, \eqref{h18}, and \eqref{h17} complete the induction step. Induction then completes the proof of \cref{c02}.
\end{proof}

\section{Perturbation lemma}
\label{b37b}

In \cref{b37} below we approximate the solution to the SFPE \eqref{b05a} through the solution to the SFPE
\eqref{b05c}, which is defined with respect to  approximating functions in the definition of the corresponding SFPE. Later in the proof of \cref{d23} these approximating functions will be the ones which can be represented by DNNs as stated in the assumption of \cref{d23}.
\begin{proposition}[Perturbation result]\label{b37}
Let $d\in \N$, 
$T, \varepsilon\in (0,\infty)$, 
$\exponentV,\exponentZ, \exponentX\in [1,\infty)$,
$q\in [1,2)$,
$c\in [1,\infty)$, $(L_i)_{i\in [0,d]\cap \Z}\in \R^{d+1}$ satisfy that $
\sum_{i=0}^{d}L_i\leq c$, $
 \frac{2}{\exponentV}+\frac{1}{\exponentX}+\frac{1}{\exponentZ}\leq 1$, $2q\leq \exponentV$,
and $\frac{1}{2q}+\frac{1}{\exponentZ}\leq 1$.
Let $\lVert\cdot \rVert\colon \R^d\to [0,\infty)$ be a norm on $\R^d$.
Let 
 $\Lambda=(\Lambda_{\nu})_{\nu\in [0,d]\cap\Z}\colon [0,T]\to \R^{1+d}$ satisfy for all $t\in [0,T]$ that $\Lambda(t)=(1,\sqrt{t},\ldots,\sqrt{t})$. Let $\pr=(\pr_\nu)_{\nu\in [0,d]\cap\Z}\colon \R^{d+1}\to\R$ satisfy
for all 
$w=(w_\nu)_{\nu\in [0,d]\cap\Z}$,
$i\in [0,d]\cap\Z$ that
$\pr_i(w)=w_i$. Let $f, \tilde{f}\in C( \R^{d+1},\R)$, $g, \tilde{g}\in C(\R^d,\R)$, 
$V\in C([0,T]\times \R^d, [0,\infty) )$ satisfy that
$\max \{ c,48e^{86c^6T^3}\}\leq V$.
Let $ (\Omega,\mathcal{F},\P)$ be a probability space.
Let $(X^{s,x}_t)_{s\in [0,T],t\in[s,T],x\in\R^d}\colon \{(\sigma,\tau)\in [0,T]^2\colon \sigma\leq \tau\}\times \R^d\times \Omega\to\R^d $,
$(Z^{s,x}_t)_{s\in [0,T),t\in(s,T],x\in\R^d}\colon \{(\sigma,\tau)\in [0,T]^2\colon \sigma< \tau\}\times \R^d\times \Omega\to\R^{d+1} $ be measurable. Assume
for all 
$i\in [0,d]\cap\Z$,
$s\in [0,T)$,
$t\in [s,T)$, $r\in (t,T]$,
 $x,y\in \R^d$, $w_1,w_2\in \R^{d+1}$,
$A\in (\mathcal{B}(\R^d))^{\otimes \R^d}$,
$B\in 
(\mathcal{B}(\R^d))^{\otimes ([t,T)\times\R^d)}$
 that
\begin{align}\label{b04b}
\max\!\left\{
\lvert g(x)\rvert, \lvert \tilde{g}(x)\rvert\right\}\leq V(T,x),\quad
\max\!\left\{
\lvert Tf(0)\rvert, \lvert T\tilde{f}(0)\rvert\right\}\leq V(t,x),
\end{align}
\begin{align}
\max\!\left\{
\lvert
f(w_1)-f(w_2)\rvert,
\lvert
\tilde{f}(w_1)-\tilde{f}(w_2)\rvert\right\}\leq \sum_{\nu=0}^{d}\left[
L_\nu\Lambda_\nu(T)
\lvert\pr_\nu(w_1-w_2) \rvert\right]
,\label{c01}
\end{align}
\begin{align}\label{b02c}
\left\lVert
V(r,X^{t,x}_r)
\right\rVert_{\exponentV}\leq V(t,x),\quad 
\left\lVert
\pr_i(Z^{t,x}_r)\right\rVert_{\exponentZ}\leq \frac{c}{\Lambda_{i}(r-t)},
\end{align}
\begin{align}
\max \!\left\{
\lvert g(x)-g(y)\rvert,
\lvert \tilde{g}(x)-\tilde{g}(y)\rvert
\right\}
\leq \frac{V(T,x)+V(T,y)}{2}
\frac{\lVert x-y\rVert}{\sqrt{T}},\label{b19}
\end{align}
\begin{align}
\left\lVert\left\lVert
X^{t,{x}}_{r}-
X^{t,{y}}_{r}\right\rVert
\right\rVert_{\exponentX}\leq 
c
\lVert x-y\rVert,\label{b21}
\end{align}
\begin{align}
\left\lVert
\pr_i\! \left(
Z^{t,{x}}_{r} 
-Z^{t,{y}}_{r}\right)\right\rVert_{\exponentZ}\leq
\frac{V(t,x)+V(t,y)}{2} \frac{\lVert x-y\rVert}{\sqrt{T}\Lambda_{i}(r-t)},\label{b39}
\end{align}
\begin{align}
\P \!\left( X_r^{t, X_t^{s,x}} = X^{s,x}_r\right)=1,\quad 
\P\!\left(X_t^{s,(\cdot)} \in A, X^{t,(\cdot)}_{(\cdot)} \in B\right) =\P\!\left (X_t^{s,(\cdot)} \in A\right)\P\!\left(X^{t,(\cdot)}_{(\cdot)} \in B\right),\label{c40}
\end{align}
\begin{align}\label{c74}
\left\lVert
\pr_i(Z^{t,x}_r-Z^{s,x}_r)\right\rVert_{\exponentZ}\leq \frac{V(t,x)+V(s,x)}{2}\frac{\sqrt{t-s}}{\sqrt{r-t}\Lambda_i(r-s)},
\end{align}
\begin{align}\label{c73}
\left\lVert
\left\lVert
X^{s,x}_t-x
\right\rVert\right\rVert_{\exponentX}\leq V(s,x)\sqrt{t-s},
\end{align}
\begin{align}
\left\lvert
\tilde{f}(w)-f(w)\right\rvert\leq \frac{\varepsilon cd^c}{T}\left(1+\sum_{\nu=0}^{d}(\Lambda_\nu(T))^q\lvert\pr_\nu(w)\rvert^q\right),\label{k03}
\end{align}
\begin{align}\label{k04}
\left\lvert
g(x)-\tilde{g}(x)\right\rvert\leq \varepsilon cd^cV(T,x),
\end{align}
and $\P(X^{s,x}_s=x)=1$.
Then the following items hold.
\begin{enumerate}[(i)]
\item \label{b17a}There exist  unique continuous functions $u, \tilde{u}\colon [0,T)\times \R^d\to \R^{d+1}$ such that for all $t\in [0,T)$, $x\in \R^d$
we have that
\begin{align}\label{c01a}
\max_{\nu\in [0,d]\cap\Z}\sup_{\tau\in [0,T), \xi\in \R^d}
\left[\Lambda_\nu(T-\tau)\frac{\lvert\pr_\nu(u(\tau,\xi))\rvert}{V(\tau ,\xi)}\right]<\infty,
\end{align}
\begin{align}\label{c01b}
\max_{\nu\in [0,d]\cap\Z}\sup_{\tau\in [0,T), \xi\in \R^d}
\left[\Lambda_\nu(T-\tau)\frac{\lvert\pr_\nu(\tilde{u}(\tau,\xi))\rvert}{V(\tau ,\xi)}\right]<\infty,
\end{align}
\begin{align}\max_{\nu\in [0,d]\cap\Z}\left[
\E\!\left [\left\lvert g(X^{t,x}_{T} )\pr_\nu(Z^{t,x}_{T})\right\rvert \right] + \int_{t}^{T}
\E \!\left[\left\lvert
f(u(r,X^{t,x}_{r}))\pr_\nu(Z^{t,x}_{r})\right\rvert\right]dr\right]<\infty,
\end{align}
\begin{align}
\max_{\nu\in [0,d]\cap\Z}\left[
\E\!\left [\left\lvert \tilde{g}(X^{t,x}_{T} )\pr_\nu(Z^{t,x}_{T})\right\rvert \right] + \int_{t}^{T}
\E \!\left[\left\lvert
\tilde{f}(\tilde{u}(r,X^{t,x}_{r}))\pr_\nu(Z^{t,x}_{r})\right\rvert\right]dr\right]<\infty,
\end{align}
\begin{align}
u(t,x)=\E\!\left [g(X^{t,x}_{T} )Z^{t,x}_{T} \right] + \int_{t}^{T}
\E \!\left[
f(u(r,X^{t,x}_{r}))Z^{t,x}_{r}\right]dr,\label{b05a}
\end{align}
and
\begin{align}
\tilde{u}(t,x)=\E\!\left [\tilde{g}(X^{t,x}_{T} )Z^{t,x}_{T} \right] + \int_{t}^{T}
\E \!\left[
\tilde{f}(\tilde{u}(r,X^{t,x}_{r}))Z^{t,x}_{r}\right]dr.\label{b05c}
\end{align}

\item For all
 $t\in [0,T)$, $x\in \R^d$ we have that
\begin{align}
\max_{i\in [0,d]\cap\Z}
\left\lVert
\Lambda_{i}(T-t)
\pr_i(\tilde{u}(t,x) 
-u(t,x))\right\rVert_2
\leq 
10 c^2d^{2c}\varepsilon V^{3q+1}(t,x)\mathrm{B}(1-\tfrac{q}{2},\tfrac{1}{2}).\label{b05d}
\end{align}
\end{enumerate}
\end{proposition}
\begin{proof}
[Proof of \cref{b37}]\cite[Lemma 2.7]{NNW2023} and the assumptions of \cref{b37} show \eqref{b17a}.

Next, \eqref{b04b}, H\"older's inequality, the fact that
$\frac{1}{\exponentV}+\frac{1}{\exponentZ}\leq 1$, \eqref{b02c}, 
Jensen's inequality,
the fact that $2q\leq \exponentV$, and the fact that 
$c\leq V $ imply for all $i\in [0,d]\cap\Z$, $s\in [0,T)$, $t\in [s,T)$, $x\in \R^d$ that
\begin{align}
&
\left\lVert
\Lambda_i(T-t)\left\lVert g(X^{t,\tilde{x}}_{T} )\pr_i(Z^{t,\tilde{x}}_{T}) \right\rVert_1\Bigr|_{\tilde{x}=X^{s,x}_t}
\right\rVert_{2q}\nonumber\\
&\leq \xeqref{b04b}
\left\lVert 
\Lambda_i(T-t)
\left\lVert V(T,X^{t,\tilde{x}}_{T} )\pr_i(Z^{t,\tilde{x}}_{T})\right\rVert_1\Bigr|_{\tilde{x}=X^{s,x}_t}
\right\rVert_{2q}\nonumber\\
&\leq 
\left\lVert 
\Lambda_i(T-t)
\left\lVert V(T,X^{t,\tilde{x}}_{T})\right\rVert_{\exponentV} \left\lVert\pr_i(Z^{t,\tilde{x}}_{T})\right\rVert_{\exponentZ}\Bigr|_{\tilde{x}=X^{s,x}_t}
\right\rVert_{2q}\nonumber\\
&\leq \xeqref{b02c}
\left\lVert\Lambda_i(T-t)V(t,\tilde{x})\frac{c}{\Lambda_i(T-t)}\Bigr|_{\tilde{x}=X^{s,x}_t} \right\rVert_{2q}\nonumber\\
&\leq c\left\lVert V(T, X^{s,x}_T)\right\rVert_{\exponentV}\nonumber\\
&\leq \xeqref{b02c} V^2(s,x).\label{k01}
\end{align}
Next,
Jensen's inequality,
\eqref{b02c}, the fact that $1\leq \exponentZ$,  \eqref{b04b}, and the fact that $c\leq V$
prove for all 
$i\in [0,d]\cap\Z$, $s\in [0,T)$, $t\in [s,T)$, $x\in \R^d$
that
\begin{align}
\int_{t}^{T}
\left\lVert
\left\lVert
\Lambda_i(T-t)
f(0)\pr_i(Z^{t,\tilde{x}}_{r})\right\rVert_1\Bigr|_{\tilde{x}=X^{s,x}_t}\right\rVert_{2q}dr
&\leq  \int_{t}^{T}\left\lvert\Lambda_i(T-t)f(0)\frac{c}{\Lambda_i(r-t)}\right\rvert dr\nonumber\\
&\leq c\lvert f(0)\rvert\int_{t}^{T}\frac{\sqrt{T-t}}{\sqrt{r-t}}\,dr\nonumber\\
&= c\lvert f(0)\rvert\sqrt{T-t}(2\sqrt{r-t}|_{r=t}^T)\nonumber\\
&\leq 2c\lvert Tf(0)\rvert\leq 2c V(s,x)\leq 2V^2(s,x).\label{k02}
\end{align}
Next, H\"older's inequality, the fact that
$\frac{1}{2q}+\frac{1}{\exponentZ}\leq 1$, the triangle inequality, \eqref{b02c}, the disintegration theorem, the flow property in \eqref{c40}, and the fact that $\sum_{\nu=0}^{d}L_\nu\leq c$
show for all 
$i\in [0,d]\cap\Z$, $s\in [0,T)$, $t\in [s,T)$, $x\in \R^d$ that
\begin{align}
&
\int_{t}^{T}\left\lVert\Lambda_i(T-t)
\sum_{\nu=0}^{d}L_\nu \Lambda_\nu (T)\left\lVert\pr_\nu(u(r,X^{t,\tilde{x}}_r))\pr_i(Z^{t,\tilde{x}}_{r})\right\rVert_1\Bigr|_{\tilde{x}=X^{s,x}_t}\right\rVert_{2q}\nonumber\\
&\leq 
\int_{t}^{T}\left\lVert\Lambda_i(T-t)
\sum_{\nu=0}^{d}L_\nu\frac{\sqrt{T}}{\sqrt{T-r}} \Lambda_\nu(T-r) \left\lVert\pr_\nu(u(r,X^{t,\tilde{x}}_r))
\right\rVert_{2q}
\left\lVert
\pr_i(Z^{t,\tilde{x}}_{r})\right\rVert_{\exponentZ}\Bigr|_{\tilde{x}=X^{s,x}_t}\right\rVert_{2q}\nonumber\\
&\leq 
\int_{t}^{T}\Lambda_i(T-t)
\sum_{\nu=0}^{d}L_\nu\frac{\sqrt{T}}{\sqrt{T-r}} \Lambda_\nu(T-r) \left\lVert\left\lVert\pr_\nu(u(r,X^{t,\tilde{x}}_r))
\right\rVert_{2q}
\left\lVert
\pr_i(Z^{t,\tilde{x}}_{r})\right\rVert_{\exponentZ}\Bigr|_{\tilde{x}=X^{s,x}_t}\right\rVert_{2q}\nonumber\\
&\leq \xeqref{b02c}
\int_{t}^{T}\Lambda_i(T-t)
\sum_{\nu=0}^{d}L_\nu\frac{\sqrt{T}}{\sqrt{T-r}} \Lambda_\nu(T-r) \left\lVert\left\lVert\pr_\nu(u(r,X^{t,\tilde{x}}_r))
\right\rVert_{2q}
\frac{c}{\Lambda_i(r-t)}\Bigr|_{\tilde{x}=X^{s,x}_t}\right\rVert_{2q}\nonumber\\
&\leq 
\int_{t}^{T}\Lambda_i(T-t)
\sum_{\nu=0}^{d}L_\nu\frac{\sqrt{T}}{\sqrt{T-r}} \Lambda_\nu(T-r)
\left\lVert \pr_\nu(u(r,X^{s,x}_r))\right\rVert_{2q}
\frac{c}{\Lambda_i(r-t)}\,dr\nonumber\\
&\leq 
\int_{t}^{T}\sqrt{T-t}
c\frac{\sqrt{T}}{\sqrt{T-r}}\max_{\nu\in [0,d]\cap\Z}\left[ \Lambda_\nu(T-r)
\left\lVert \pr_\nu(u(r,X^{s,x}_r))\right\rVert_{2q}\right]
\frac{c}{\sqrt{r-t}}\,dr\nonumber\\
&\leq \int_{t}^{T}c^2T
\max_{\nu\in [0,d]\cap\Z}\left[ \Lambda_\nu(T-r)
\left\lVert \pr_\nu(u(r,X^{s,x}_r))\right\rVert_{2q}\right]\frac{dr}{\sqrt{(T-r)(r-t)}}.
\end{align}
This,
\eqref{b05a}, the triangle inequality, \eqref{c01}, \eqref{k01}, and \eqref{k02} imply for all $i\in [0,d]\cap\Z$, $s\in [0,T)$, $t\in [s,T)$, $x\in \R^d$ that
\begin{align}
&\Lambda_i(T-t) \left\lVert\pr_i(u(t,X_t^{s,x}))\right\rVert_{2q}\nonumber\\
&
=\left\lVert\Lambda_i(T-t)\E\!\left [g(X^{t,\tilde{x}}_{T} )\pr_i(Z^{t,\tilde{x}}_{T}) \right]\Bigr|_{\tilde{x}=X^{s,x}_t}+
\Lambda_i(T-t)
\int_{t}^{T}
\E \!\left[
f(u(r,X^{t,\tilde{x}}_{r}))\pr_i(Z^{t,\tilde{x}}_{r})\right]\Bigr|_{\tilde{x}=X^{s,x}_t}dr\right\rVert_{2q}\nonumber\\
&\leq 
\left\lVert
\Lambda_i(T-t)\left\lVert g(X^{t,\tilde{x}}_{T} )\pr_i(Z^{t,\tilde{x}}_{T}) \right\rVert_1\Bigr|_{\tilde{x}=X^{s,x}_t}
\right\rVert_{2q}
+
\int_{t}^{T}
\left\lVert
\left\lVert
\Lambda_i(T-t)
f(u(r,X^{t,\tilde{x}}_{r}))\pr_i(Z^{t,\tilde{x}}_{r})\right\rVert_1\Bigr|_{\tilde{x}=X^{s,x}_t}\right\rVert_{2q}dr\nonumber\\
&\leq 
\left\lVert
\Lambda_i(T-t)\left\lVert g(X^{t,\tilde{x}}_{T} )\pr_i(Z^{t,\tilde{x}}_{T}) \right\rVert_1\Bigr|_{\tilde{x}=X^{s,x}_t}
\right\rVert_{2q}
+
\int_{t}^{T}
\left\lVert
\left\lVert
\Lambda_i(T-t)
f(0)\pr_i(Z^{t,\tilde{x}}_{r})\right\rVert_1\Bigr|_{\tilde{x}=X^{s,x}_t}\right\rVert_{2q}dr\nonumber\\
&\quad +\int_{t}^{T}\left\lVert\Lambda_i(T-t)
\sum_{\nu=0}^{d}L_\nu \Lambda_\nu (T)\left\lVert\pr_\nu(u(r,X^{t,\tilde{x}}_r))\pr_i(Z^{t,\tilde{x}}_{r})\right\rVert_1\Bigr|_{\tilde{x}=X^{s,x}_t}\right\rVert_{2q}\nonumber\\
&\leq \xeqref{k01} V^2(s,x)+\xeqref{k02}2V^2(s,x)+ \int_{t}^{T}c^2T
\max_{\nu\in [0,d]\cap\Z}\left[ \Lambda_\nu(T-r)
\left\lVert \pr_\nu(u(r,X^{s,x}_r))\right\rVert_{2q}\right]\frac{dr}{\sqrt{(T-r)(r-t)}}\nonumber\\
&\leq 3V^2(s,x)
+ \int_{t}^{T}c^2T
\max_{\nu\in [0,d]\cap\Z}\left[ \Lambda_\nu(T-r)
\left\lVert \pr_\nu(u(r,X^{s,x}_r))\right\rVert_{2q}\right]\frac{dr}{\sqrt{(T-r)(r-t)}}.
\end{align}
This proves for all $s\in [0,T)$, $t\in [s,T)$, $x\in \R^d$ that
\begin{align}
&\max_{\nu\in [0,d]\cap\Z}\left[ \Lambda_\nu(T-t)
\left\lVert \pr_\nu(u(t,X^{s,x}_t))\right\rVert_{2q}\right]\nonumber\\
&\leq 3V^2(s,x)
+ \int_{t}^{T}c^2T
\max_{\nu\in [0,d]\cap\Z}\left[ \Lambda_\nu(T-r)
\left\lVert \pr_\nu(u(r,X^{s,x}_r))\right\rVert_{2q}\right]\frac{dr}{\sqrt{(T-r)(r-t)}}.\label{k05b}
\end{align}
This, \eqref{c01a}, \eqref{b02c}, the fact that $2q\leq \exponentV $, and Gr\"onwall's inequality (see  \cite[Corollary~2.4]{NNW2023}) show for all $ s\in [0,T)$, $t\in [s,T)$, $x\in \R^d$ that
\begin{align}
\max_{\nu\in [0,d]\cap\Z}\left[ \Lambda_\nu(T-t)
\left\lVert \pr_\nu(u(t,X^{s,x}_t))\right\rVert_{2q}\right]\leq 
6V^2(s,x)e^{86c^6T^3}.\label{k05}
\end{align}
Next, \eqref{k04}, H\"older's inequality, the fact that $\frac{1}{\exponentV}+\frac{1}{\exponentZ}\leq 1$, 
\eqref{b02c},
Jensen's inequality, and the fact that $\frac{1}{\exponentV}\leq \frac{1}{2}$
imply for all 
$s\in [0,T)$, $t\in [s,T)$, $x\in \R^d$, $\nu\in [0,d]\cap\Z $ that
\begin{align}
& 
\Lambda_\nu(T-t)
\left\lVert
 \left\lVert
(\tilde{g}(X^{t,\tilde{x}}_T)-g(X^{t,\tilde{x}}_T))\pr_\nu(Z^{t,\tilde{x}}_T)
\right\rVert_{1}\Bigr|_{\tilde{x}=X^{s,x}_t}
\right\rVert_{2}\nonumber\\
&\leq 
\Lambda_\nu(T-t)
\left\lVert
 \left\lVert
\varepsilon cd^c V(T, X^{t,\tilde{x}}_T) \pr_\nu(Z^{t,\tilde{x}}_T)
\right\rVert_{1}\Bigr|_{\tilde{x}=X^{s,x}_t}
\right\rVert_{2}\nonumber\\
&\leq \Lambda_\nu(T-t)
\left\lVert
\varepsilon cd^c 
\left\lVert
V(T, X^{t,\tilde{x}}_T)
\right\rVert_{\exponentV}
\left\lVert
 \pr_\nu(Z^{t,\tilde{x}}_T)\right\rVert_{\exponentZ}
\Bigr|_{\tilde{x}=X^{s,x}_t}
\right\rVert_{2}\nonumber\\
&\leq 
\Lambda_\nu(T-t)
\left\lVert
\varepsilon cd^c 
V(t,X_t^{s,x})
\frac{c}{\Lambda_\nu(T-t)}
\right\rVert_{\exponentV}\nonumber\\
&\leq \varepsilon c^2 d^cV(s,x).\label{h01}
\end{align}
Next, \eqref{c01}, H\"older's inequality, the fact that
$\frac{1}{2}+\frac{1}{\exponentZ}\leq 1$, \eqref{b02c}, the disintegration theorem, \eqref{c40}, and the fact that
$
\sum_{i=0}^{d}L_i\leq c$ prove for all $\nu\in[0,d]\cap\Z$,
$ s\in [0,T)$, $t\in [s,T)$, $x\in \R^d$ that
\begin{align}
&
\int_{t}^{T}
\Lambda_\nu(T-t)
\left\lVert
\left\lVert\tilde{f}(\tilde{u}(r,X^{t,\tilde{x}}_r))\pr_\nu(Z^{t,\tilde{x}}_{r})-
\tilde{f}({u}(r,X^{t,\tilde{x}}_r))\pr_\nu(Z^{t,\tilde{x}}_{r})\right\rVert_1\Bigr|_{\tilde{x}=X^{s,x}_t}\right\rVert_2
dr\nonumber\\
&\leq \xeqref{c01}
\int_{t}^{T}
\Lambda_\nu(T-t)
\left\lVert
\left\lVert
\sum_{i=0}^{d}L_i\Lambda_{i}(T)\left\lvert\pr_i(\tilde{u}(r,X^{t,\tilde{x}}_r) 
-u(r,X^{t,\tilde{x}}_r))\right\rvert\left\lvert
\pr_\nu(Z^{t,\tilde{x}}_{r})\right\rvert\right\rVert_1\Biggr|_{\tilde{x}=X^{s,x}_t}\right\rVert_2
dr\nonumber\\
&\leq 
\int_{t}^{T}
\Lambda_\nu(T-t)
\sum_{i=0}^{d}L_i\frac{\sqrt{T}}{\sqrt{T-r}}
\left\lVert
\left\lVert
\Lambda_{i}(T-r)
\pr_i(\tilde{u}(r,X^{t,\tilde{x}}_r) 
-u(r,X^{t,\tilde{x}}_r))
\pr_\nu(Z^{t,\tilde{x}}_{r})\right\rVert_1\Bigr|_{\tilde{x}=X^{s,x}_t}\right\rVert_2
dr\nonumber\\
&\leq 
\int_{t}^{T}
\Lambda_\nu(T-t)
\sum_{i=0}^{d}L_i\frac{\sqrt{T}}{\sqrt{T-r}}
\left\lVert
\left\lVert
\Lambda_{i}(T-r)
\pr_i(\tilde{u}(r,X^{t,\tilde{x}}_r) 
-u(r,X^{t,\tilde{x}}_r))\right\rVert_2
\left\lVert
\pr_\nu(Z^{t,\tilde{x}}_{r})\right\rVert_{\exponentZ}\Bigr|_{\tilde{x}=X^{s,x}_t}\right\rVert_2
dr\nonumber\\
&\leq
\int_{t}^{T}
\Lambda_\nu(T-t)
\sum_{i=0}^{d}L_i\frac{\sqrt{T}}{\sqrt{T-r}}
\left\lVert
\left\lVert
\Lambda_{i}(T-r)
\pr_i(\tilde{u}(r,X^{t,\tilde{x}}_r) 
-u(r,X^{t,\tilde{x}}_r))\right\rVert_2
 \xeqref{b02c}
\frac{c}{\Lambda_\nu(r-t)}\Bigr|_{\tilde{x}=X^{s,x}_t}\right\rVert_2
dr\nonumber\\
&\leq 
\int_{t}^{T}\sqrt{T-t}
c\frac{\sqrt{T}}{\sqrt{T-r}}
\max_{i\in [0,d]\cap\Z}
\left\lVert
\Lambda_{i}(T-r)
\pr_i(\tilde{u}(r,X^{s,x}_r) 
-u(r,X^{s,x}_r))\right\rVert_2\frac{c}{\sqrt{r-t}}\,
dr\nonumber\\
&\leq 
\int_{t}^{T}c^2 T
\max_{i\in [0,d]\cap\Z}
\left\lVert
\Lambda_{i}(T-r)
\pr_i(\tilde{u}(r,X^{s,x}_r) 
-u(r,X^{s,x}_r))\right\rVert_2\frac{dr}{\sqrt{(T-r)(r-t)}}.
\label{h02}
\end{align}
Next,
\eqref{b02c} shows for all $\nu\in[0,d]\cap\Z$,
$ s\in [0,T)$, $t\in [s,T)$, $x\in \R^d$ 
that
\begin{align}
\int_{t}^{T}
\Lambda_\nu(T-t)
\left\lVert
\left\lVert
 \varepsilon
\pr_\nu(Z^{t,\tilde{x}}_{r})
\right\rVert_1\bigr|_{\tilde{x}=X^{s,x}_t}\right\rVert_2
dr&\leq \int_{t}^{T}\Lambda_{\nu}(T-t)\varepsilon \frac{c}{\Lambda_{\nu}(r-t)}dr
\nonumber\\
&\leq \int_{t}^{T}\varepsilon c\frac{\sqrt{T-t}}{\sqrt{r-t}}\,dr\nonumber\\
&=\varepsilon c\sqrt{T-t}(2\sqrt{r-t}|_{r=t}^T)\nonumber\\
&\leq 2\varepsilon c T.\label{k06}
\end{align}
Next, the substitution 
$s=\frac{r-t}{T-t}$, $ds=\frac{dr}{T-t}$ implies for all $t\in [0,T)$ that
\begin{align}
\int_{t}^{T}\frac{(T-t)^\frac{1}{2}T^\frac{q}{2}}{(T-r)^\frac{q}{2}(r-t)^\frac{1}{2}}\,dr
&=
\int_{0}^{1}\frac{(T-t)^\frac{1}{2}T^\frac{q}{2}(T-t)}{[(T-t)(1-s)]^\frac{q}{2}[(T-t)s]^\frac{1}{2}}\,ds\nonumber\\
&
\leq T\int_{0}^{1}\frac{ds}{(1-s)^\frac{q}{2}s^\frac{1}{2}}
=T\mathrm{B}(1-\tfrac{q}{2}, \tfrac{1}{2}).
\end{align}
This, H\"older's inequality, the fact that 
$\frac{1}{2}+\frac{1}{\exponentZ}\leq 1$,
\eqref{b02c}, the disintegration theorem, the flow property in \eqref{c40}, \eqref{k05}, and the fact that $ 6e^{86c^6T^3}\leq V$ prove for all 
$\nu\in[0,d]\cap\Z$,
$ s\in [0,T)$, $t\in [s,T)$, $x\in \R^d$ that
\begin{align}
&
\int_{t}^{T}
\Lambda_\nu(T-t)
\left\lVert
\left\lVert
 \varepsilon\left(\sum_{i=0}^{d}(\Lambda_i(T))^q(\pr_i({u}(r,X^{t,\tilde{x}}_r)))^q\right)
\pr_\nu(Z^{t,\tilde{x}}_{r})
\right\rVert_1\Biggr|_{\tilde{x}=X^{s,x}_t}\right\rVert_2
dr\nonumber\\
&\leq 
\int_{t}^{T}
\Lambda_\nu(T-t)\varepsilon
\sum_{i=0}^{d}
\frac{{T}^\frac{q}{2}}{(T-r)^\frac{q}{2}}
(\Lambda_i(T-r))^q
\left\lVert
 \left\lVert
(\pr_i({u}(r,X^{t,\tilde{x}}_r)))^q
\pr_\nu(Z^{t,\tilde{x}}_{r})
\right\rVert_1\bigr|_{\tilde{x}=X^{s,x}_t}\right\rVert_2
dr\nonumber\\
&\leq 
\int_{t}^{T}
\Lambda_\nu(T-t)\varepsilon
\sum_{i=0}^{d}
\frac{{T}^\frac{q}{2}}{(T-r)^\frac{q}{2}}
(\Lambda_i(T-r))^q
\left\lVert
 \left\lVert
(\pr_i({u}(r,X^{t,\tilde{x}}_r)))^q
\right\rVert_2
\left\lVert\pr_\nu(Z^{t,\tilde{x}}_{r})
\right\rVert_\exponentZ\bigr|_{\tilde{x}=X^{s,x}_t}\right\rVert_2
dr\nonumber\\
&\leq 
\int_{t}^{T}
\Lambda_\nu(T-t)\varepsilon
\sum_{i=0}^{d}
\frac{{T}^\frac{q}{2}}{(T-r)^\frac{q}{2}}
(\Lambda_i(T-r))^q
\left\lVert
 \left\lVert
(\pr_i({u}(r,X^{t,\tilde{x}}_r)))^q
\right\rVert_2\xeqref{b02c}
\frac{c}{\Lambda_{\nu}(T-t)}\Bigr|_{\tilde{x}=X^{s,x}_t}\right\rVert_2
dr\nonumber\\
&\leq \int_{t}^{T}
\Lambda_\nu(T-t)\varepsilon
\sum_{i=0}^{d}
\frac{{T}^\frac{q}{2}}{(T-r)^\frac{q}{2}}
(\Lambda_i(T-r))^q
\xeqref{c40}
 \left\lVert \pr_i({u}(r,X^{s,x}_r))
\right\rVert^q_{2q}
\frac{c}{\Lambda_{\nu}(r-t)}\,
dr\nonumber\\
&\leq 
\int_{t}^{T}
\sqrt{T-t}\varepsilon
(d+1)
\frac{{T}^\frac{q}{2}}{(T-r)^\frac{q}{2}}
\left[\xeqref{k05}6V^2(s,x)e^{86c^6T^3}\right]^q
\frac{c}{\sqrt{r-t}}\,
dr\nonumber\\
&\leq c\varepsilon(d+1) V^{3q}(s,x)\int_{t}^{T}\frac{(T-t)^\frac{1}{2}T^\frac{q}{2}}{(T-r)^\frac{q}{2}(r-t)^\frac{1}{2}}\,dr\nonumber\\
&\leq c\varepsilon(d+1) V^{3q}(s,x)T\mathrm{B}(1-\tfrac{q}{2},\tfrac{1}{2}).
\end{align}
This, \eqref{k03}, the triangle inequality,  \eqref{k06}, and the fact that
$\mathrm{B}(1-\frac{q}{2}, \frac{1}{2})\geq 1 $
 show for
all $\nu\in[0,d]\cap\Z$,
$ s\in [0,T)$, $t\in [s,T)$, $x\in \R^d$  that
\begin{align}
&
\int_{t}^{T}
\Lambda_\nu(T-t)
\left\lVert
\left\lVert
\tilde{f}({u}(r,X^{t,\tilde{x}}_r))\pr_\nu(Z^{t,\tilde{x}}_{r})
-{f}({u}(r,X^{t,\tilde{x}}_r))\pr_\nu(Z^{t,\tilde{x}}_{r})
\right\rVert_1\Bigr|_{\tilde{x}=X^{s,x}_t}\right\rVert_2
dr\nonumber\\
&\leq 
\int_{t}^{T}
\Lambda_\nu(T-t)
\left\lVert
\left\lVert
 \frac{\varepsilon cd^c}{T}\left(1+\sum_{\nu=0}^{d}(\Lambda_i(T))^q(\pr_\nu({u}(r,X^{t,\tilde{x}}_r)))^q\right)
\pr_\nu(Z^{t,\tilde{x}}_{r})
\right\rVert_1\Biggr|_{\tilde{x}=X^{s,x}_t}\right\rVert_2
dr\nonumber\\
& 
\leq \frac{cd^c}{T} \int_{t}^{T}
\Lambda_\nu(T-t)
\left\lVert
\left\lVert
 \varepsilon 
\pr_\nu(Z^{t,\tilde{x}}_{r})
\right\rVert_1\Bigr|_{\tilde{x}=X^{s,x}_t}\right\rVert_2
dr\nonumber\\
&\quad +\frac{cd^c}{T}\int_{t}^{T}
\Lambda_\nu(T-t)
\left\lVert
\left\lVert
 \varepsilon\left(\sum_{i=0}^{d}(\Lambda_\nu(T))^q(\pr_i({u}(r,X^{t,\tilde{x}}_r)))^q\right)
\pr_\nu(Z^{t,\tilde{x}}_{r})
\right\rVert_1\Biggr|_{\tilde{x}=X^{s,x}_t}\right\rVert_2
dr\nonumber\\
&\leq \frac{cd^c}{T} \left(2\varepsilon c T+c\varepsilon(d+1) V^{3q}(s,x)T\mathrm{B}(1-\tfrac{q}{2},\tfrac{1}{2})\right)
\nonumber\\
&\leq cd^cc\varepsilon(d+3) V^{3q}(s,x)\mathrm{B}(1-\tfrac{q}{2},\tfrac{1}{2})\nonumber\\
&\leq 4dc^2d^c\varepsilon V^{3q}(s,x)\mathrm{B}(1-\tfrac{q}{2},\tfrac{1}{2})\nonumber\\
&\leq 4c^2d^{2c}\varepsilon V^{3q}(s,x)\mathrm{B}(1-\tfrac{q}{2},\tfrac{1}{2}).
\end{align}
This, \eqref{b05a}, \eqref{b05c}, the triangle inequality, \eqref{h01}, \eqref{h02}, the fact that $V\leq V^{3q}$, and the fact that
$\mathrm{B}(1-\tfrac{q}{2},\tfrac{1}{2})\geq 1$ prove for all
$\nu\in [0,d]\cap\Z$, 
$ s\in [0,T)$, $t\in [s,T)$, $x\in \R^d$  that
\begin{align}
&
\Lambda_\nu(T-t)
\left\lVert\pr_\nu\!\left(
u(t,X^{s,x}_t)-\tilde{u}(t,X^{s,x}_t)\right)\right\rVert_2\nonumber\\
&
\leq \Lambda_\nu(T-t)
\left\lVert \E \!\left[\tilde{g}(X^{t,\tilde{x}}_T)\pr_\nu(Z_T^{t,\tilde{x}})
-{g}(X^{t,\tilde{x}}_T)\pr_\nu(Z_T^{t,\tilde{x}})
\right]\Bigr|_{\tilde{x}= X^{s,x}_t}
\right\rVert_2\nonumber\\
&\quad +\int_{t}^{T}
\Lambda_\nu(T-t)
\left\lVert
\left\lvert\E \!\left[\tilde{f}(\tilde{u}(r,X^{t,\tilde{x}}_r))\pr_\nu(Z^{t,\tilde{x}}_{r})\right]-
\E \!\left[{f}({u}(r,X^{t,\tilde{x}}_r))\pr_\nu(Z^{t,\tilde{x}}_{r})\right]\right\rvert\Bigr|_{\tilde{x}=X^{s,x}_t}
\right\rVert_2
dr\nonumber\\
&\leq 
\Lambda_\nu(T-t)
\left\lVert
 \left\lVert
(\tilde{g}(X^{t,\tilde{x}}_T)-g(X^{t,\tilde{x}}_T))\pr_\nu(Z^{t,\tilde{x}}_T)
\right\rVert_{1}\Bigr|_{\tilde{x}=X^{s,x}_t}
\right\rVert_{2}\nonumber\\
&\quad +\int_{t}^{T}
\Lambda_\nu(T-t)
\left\lVert
\left\lVert\tilde{f}(\tilde{u}(r,X^{t,\tilde{x}}_r))\pr_\nu(Z^{t,\tilde{x}}_{r})-
\tilde{f}({u}(r,X^{t,\tilde{x}}_r))\pr_\nu(Z^{t,\tilde{x}}_{r})\right\rVert_1\Bigr|_{\tilde{x}=X^{s,x}_t}\right\rVert_2
dr\nonumber\\
&\quad+
\int_{t}^{T}
\Lambda_\nu(T-t)
\left\lVert
\left\lVert
\tilde{f}({u}(r,X^{t,\tilde{x}}_r))\pr_\nu(Z^{t,\tilde{x}}_{r})
-{f}({u}(r,X^{t,\tilde{x}}_r))\pr_\nu(Z^{t,\tilde{x}}_{r})
\right\rVert_1\Bigr|_{\tilde{x}=X^{s,x}_t}\right\rVert_2
dr\nonumber\\
&\leq \xeqref{h01}
\varepsilon c^2 d^cV(s,x)+\xeqref{h02}
\int_{t}^{T}c^2 T
\max_{i\in [0,d]\cap\Z}
\left\lVert
\Lambda_{i}(T-r)
\pr_i(\tilde{u}(r,X^{s,x}_r) 
-u(r,X^{s,x}_r))\right\rVert_2\frac{dr}{\sqrt{(T-r)(r-t)}}
\nonumber
\\
&\quad +4c^2d^{2c}\varepsilon V^{3q}(s,x)\mathrm{B}(1-\tfrac{q}{2},\tfrac{1}{2})\nonumber\\
&\leq 5 c^2d^{2c}\varepsilon V^{3q}(s,x)\mathrm{B}(1-\tfrac{q}{2},\tfrac{1}{2})\nonumber\\
&\quad +
\int_{t}^{T}c^2 T
\max_{i\in [0,d]\cap\Z}
\left\lVert
\Lambda_{i}(T-r)
\pr_i(\tilde{u}(r,X^{s,x}_r) 
-u(r,X^{s,x}_r))\right\rVert_2\frac{dr}{\sqrt{(T-r)(r-t)}}.
\end{align}
Hence, 
we have for all
$ s\in [0,T)$, $t\in [s,T)$, $x\in \R^d$  that
\begin{align}
&\max_{i\in [0,d]\cap\Z}
\left\lVert
\Lambda_{i}(T-t)
\pr_i(\tilde{u}(t,X^{s,x}_t) 
-u(t,X^{s,x}_t))\right\rVert_2\nonumber\\
&\leq 5 c^2d^{2c}\varepsilon V^{3q}(s,x)\mathrm{B}(1-\tfrac{q}{2},\tfrac{1}{2})\nonumber\\
&\quad +
\int_{t}^{T}c^2 T
\max_{i\in [0,d]\cap\Z}
\left\lVert
\Lambda_{i}(T-r)
\pr_i(\tilde{u}(r,X^{s,x}_r) 
-u(r,X^{s,x}_r))\right\rVert_2\frac{dr}{\sqrt{(T-r)(r-t)}}.\label{k05c}
\end{align}
This, Gr\"onwall's inequality (see  \cite[Corollary~2.4]{NNW2023}), and the fact that 
$10e^{86c^6T^3}\leq V$ 
 imply for all
$ s\in [0,T)$, $t\in [s,T)$, $x\in \R^d$  that
\begin{align}
\max_{i\in [0,d]\cap\Z}
\left\lVert
\Lambda_{i}(T-t)
\pr_i(\tilde{u}(t,X^{s,x}_t) 
-u(t,X^{s,x}_t))\right\rVert_2
&
\leq 
10 c^2d^{2c}\varepsilon V^{3q}(s,x)\mathrm{B}(1-\tfrac{q}{2},\tfrac{1}{2})
e^{86c^6T^3}\nonumber\\
&\leq 
10 c^2d^{2c}\varepsilon V^{3q+1}(s,x)\mathrm{B}(1-\tfrac{q}{2},\tfrac{1}{2}).
\end{align}
This and the fact that
$\forall\, s\in [0,T),x\in \R^d\colon \P(X^{s,x}_s=x)=1$
complete the proof of \cref{b37}.
\end{proof}
\section{Proof of the main theorem}\label{s05}
\begin{proof}[Proof of \cref{d23}]First of all, we need to introduce an MLP setting.
Let $\varrho \colon[0,T]^2\to\R$
satisfy for all $t\in [0,T)$, $s\in (t,T)$ that
\begin{align}
\varrho(t,s)=\frac{1}{\mathrm{B}(\tfrac{1}{2},\tfrac{1}{2})}\frac{1}{\sqrt{(T-s)(s-t)}}.
\end{align}
Let $ ( \Omega, \mathcal{F}, \P )$
be a probability space.
For every $d\in \N$ let
$
  W^{d, \theta }\colon [0,T] \times \Omega \to \R^d 
$, 
$ \theta \in \Theta $,
be standard
$(\mathbb{F}_t)_{t\in[0,T]}$-Brownian motions
with continuous sample paths.
Let $\mathfrak{t}^\theta\colon \Omega\to(0,1)$, $\theta\in \Theta$, be 
i.i.d.\ random variables. Assume
for all $b\in (0,1)$
that 
\begin{align}
\P(\mathfrak{t}^0\le b)=\frac{1}{\mathrm{B}(\tfrac{1}{2},\tfrac{1}{2})}\int_0^b \frac{dr}{\sqrt{r(1-r)}}.
\end{align} 
Assume that
$(W^{d,\theta})_{d\in \N,\theta \in \Theta}$ and $(\mathfrak{t}^\theta)_{\theta \in \Theta}$
are independent. Let $\mathfrak{T}^{\theta}\colon [0,T)\times \Omega\to [0,T)$, $\theta\in\Theta$,
 satisfy for all $n \in \N_0$, $t\in [0,T)$, $\varepsilon\in (0,1)$
that $\mathfrak{T}^{\theta} _t = t+ (T-t)\mathfrak{t}^{\theta}$.
Let
$ 
  {U}_{ n,m, \varepsilon}^{d,\theta }
  \colon[0,T)\times\R^d\times\Omega\to\R^{1+d}
$,
$d\in \N$,
$n,m\in\Z$, $\theta\in\Theta$, satisfy
for all 
$
  d,n,m \in \N
$,
$ \theta \in \Theta $,
$ t\in [0,T)$,
$x \in \R^d$, $\varepsilon\in (0,1)$
that $
{U}_{-1,m,\varepsilon}^{d,\theta}(t,x)={U}_{0,m,\varepsilon}^{d,\theta}(t,x)=0$ and
\begin{equation}  \begin{split}
&{U}_{n,m,\varepsilon}^{d,\theta}(t,x)
=(g^d_\varepsilon(x), 0  ) + \sum_{i=1}^{m^n}
\left(1,\tfrac{W^{d,(\theta, 0, -i)}_{T}- W^{d,(\theta, 0, -i)}_{t}
  }{ T - t }
  \right)^\top
\tfrac{\big(g_\varepsilon^d(x+W^{d,(\theta,0,-i)}_T-W^{d,(\theta,0,-i)}_t)-g_\varepsilon^d(x)\big)}{m^n}
\\
&+\sum_{\ell=0}^{n-1}\sum_{i=1}^{m^{n-\ell}}
 \left(
  1 ,
  \tfrac{ 
  W_{\mathfrak{T}^{(\theta, \ell,i)}_t}^{d,(\theta,\ell,i)}- W^{d,(\theta,\ell, i)}_{t}
  }{ \mathfrak{T}^{(\theta, \ell,i)}_t-t}
  \right)^\top
\tfrac{  \left(f^d_\varepsilon\circ {U}_{\ell,m,\varepsilon}^{d,(\theta,\ell,i)}-\1_{\N}(\ell)f^d_\varepsilon\circ  {U}_{\ell-1,m,\varepsilon}^{d,(\theta,-\ell,i)}\right)
(\mathfrak{T}^{(\theta,\ell,i)}_t,x+W_{\mathfrak{T}^{(\theta, \ell,i)}_t}^{d,(\theta,\ell,i)}-W_t^{d,(\theta,\ell,i)})} {m^{n-\ell}\varrho(t,\mathfrak{T}^{(\theta, \ell,i)}_t)}
.
\end{split} \label{d02b} \end{equation}
For every $d\in \N$, $x\in \R^d$, $s\in [0,T]$, $t\in [s,T]$ let 
\begin{align}
a(d)=4+48e^{86c^6T^3}+8\beta d+ d^{2c},\quad \varphi_d(x)=(a(d)+\lVert x\rVert^2)^{8\beta},\label{d02}
\end{align}
\begin{align}\label{d03}
V_d(t,x)= c e^{\beta (T-t)} (\varphi_d(x))^\frac{1}{8},
\end{align}
\begin{align}\label{d07}
X_t^{d,\theta,s,x}= x+W^{d,\theta}_t-W^{d,\theta}_s.
\end{align}
For every $d\in \N$, $x\in \R^d$, $s\in [0,T)$, $t\in (s,T]$ let
\begin{align}
Z^{d,\theta,s,x}_t= \left(1,\frac{W^{d,\theta}_t-W^{d,\theta}_s}{t-s}\right)^\top .
\label{d01}
\end{align}
The fact that $8\beta d \leq \varphi^{\frac{1}{8\beta}}$
and  \cite[Lemma 2.6]{HN2022a}(applied for every $d\in \N$ with $m\gets d $, $d\gets d$, $p\gets 8\beta$, $a\gets a(d)$, $c\gets 1$, $\mu\gets 0$, $\sigma \gets \mathrm{Id}_{\R^{d\times d}}$ in the notation of \cite[Lemma 2.6]{HN2022a})  prove for all $d\in \N$ that
$\frac{1}{2}\sum_{k=1}^d\frac{\partial^2 \varphi_d}{\partial x_k^2}\leq 16\beta \varphi_d $
and hence
$\frac{1}{2}\sum_{k=1}^d\frac{\partial^2 V_d^8}{\partial x_k^2}\leq 16\beta V_d^8 $. This and the fact that
$\forall\,d\in \N\colon \frac{\partial V_d^8}{\partial t}=-8\beta V_d^8$ show for all $d\in \N$ that
\begin{align}
\frac{\partial V_d^8}{\partial t}+
\frac{1}{2}\sum_{k=1}^d\frac{\partial^2 V_d^8}{\partial x_k^2}
\leq -8\beta V_d^8+16\beta V_d^8=8\beta V_d^8.
\end{align}
This and \cite[Lemma 2.2]{CHJ2021} (applied 
for every $s\in [0,T]$, $t\in [s,T]$
with 
$V\gets V_d^8$,
$\alpha\gets 8\beta $, 
$\tau\gets T$,
$\mu\gets 0$, $\sigma\gets\mathrm{Id}_{\R^{d\times d}}$, 
$t\gets t-s $
in the notation of \cite[Lemma 2.2]{CHJ2021}) imply for all $d\in \N$, $s\in [0,T]$, $t\in [s,T]$ that
\begin{align}
e^{8\beta(T-t+s)}
\E \!\left[\varphi_d (X^{d,0,s,x}_t )\right]
&
=
\E\!\left[
e^{8\beta (T-(t-s))}\varphi_d(x+W^{d,0}_t-W^{d,0}_s)\right]\nonumber\\
&=
\E\!\left[
e^{8\beta (T-(t-s))}\varphi_d(x+W^{d,0}_{t-s})\right]\nonumber\\
&=c^{-8}
\E \!\left[
V_d^8(t-s,x+W^{d,0}_{t-s})\right]\nonumber\\&
\leq c^{-8}V^8_d(0,x)
=e^{8\beta T}\varphi_d(x),
\end{align}
hence
\begin{align}
\E \!\left[V_d^8(t,X_t^{d,0,s,x})\right]
=
\E \!\left[c^8e^{8\beta(T-t)}\varphi_d (X^{d,0,s,x}_t )\right]\leq c^8e^{8\beta(T-s)}\varphi_d(x)=V^8_d(s,x),
\end{align}
and hence
\begin{align}\label{d05}
\left\lVert
V_d(t,X_t^{d,0,s,x})\right\rVert_8\leq V_d(s,x).
\end{align}
Next, recalling the fact that $c\geq 2$ proves for all $d\in \N$, $i\in [1,d]\cap\Z$, $t\in [0,T)$, $r\in (0,T]$ that
\begin{align}
\left\lVert
\frac{W^{d,0,i}_r-W^{d,0,i}_t}{r-t}\right\rVert_8
=\frac{1}{\sqrt{r-t}}
\left\lVert
\frac{W^{d,0,i}_r-W^{d,0,i}_t}{\sqrt{r-t}}\right\rVert_8\leq \frac{c}{\sqrt{r-t}},
\end{align}
where $W^{d,0,i}$ is the $i$-th coordinate of the Brownian motion 
$W^{d,0}$.
This and \eqref{d01}  demonstrate for all $d\in \N$, 
$x\in \R^d$,
$i\in [0,d]\cap\Z$, $t\in [0,T)$, $r\in (0,T]$ that
\begin{align}
\left\lVert
\pr^d_i(Z^{d,0,t,x}_r)\right\rVert_8\leq \frac{c}{\Lambda^d_i(r-t)}.\label{d06}
\end{align}
Next, the triangle inequality implies for all 
$d\in \N$, $i\in [1,d]\cap \N$, $s\in [0,T)$, $t\in [s,T)$, $r\in (s,T]$ that
\begin{align}
&\left\lVert
\frac{W^{d,0,i}_r-W^{d,0,i}_t}{r-t}
-\frac{W^{d,0,i}_r-W^{d,0,i}_s}{r-s}\right\rVert_8\nonumber\\
&
=\left\lVert\left(\frac{1}{r-t}-\frac{1}{r-s}\right)(W^{d,0,i}_r-W^{d,0,i}_t)
+\frac{1}{r-s}\left(W^{d,0,i}_r-W^{d,0,i}_t - \left(W^{d,0,i}_r-W^{d,0,i}_s\right) \right)\right\rVert_8\nonumber\\
&\leq \frac{t-s}{(r-t)(r-s)}\left\lVert
W^{d,0,i}_r-W^{d,0,i}_t
\right\rVert_8
+\frac{1}{r-s}\left\lVert W^{d,0,i}_t -W^{d,0,i}_s \right\rVert_8\nonumber\\
&\leq \frac{t-s}{(r-t)(r-s)}105^\frac{1}{8}\sqrt{r-t}
+\frac{1}{r-s}105^\frac{1}{8}\sqrt{t-s}\nonumber\\
&\leq \frac{105^\frac{1}{8}\sqrt{t-s}}{\sqrt{r-t}\sqrt{r-s}}
+
\frac{105^\frac{1}{8}\sqrt{t-s}}{\sqrt{r-t}\sqrt{r-s}}\nonumber\\
&\leq \frac{4\sqrt{t-s}}{\sqrt{r-t}\sqrt{r-s}}.
\end{align}
This, \eqref{d01}, and the fact that $\forall\, d\in \N\colon 4\leq V_d $ show for all $d\in \N$,
$s\in [0,T)$, $t\in [s,T)$, $r\in (s,T]$ that
\begin{align}\label{c74d}
\left\lVert
\pr_i^d(Z^{d,0,t,x}_r-Z^{d,0,s,x}_r)\right\rVert_{8}\leq \frac{V_d(t,x)+V_d(s,x)}{2}\frac{\sqrt{t-s}}{\sqrt{r-t}\Lambda_i^d(r-s)}.
\end{align}
Next, \eqref{b04bb}--\eqref{k04b}, \eqref{d02}, and \eqref{d03} imply for all
$d\in \N$, $\varepsilon\in (0,1)$, $x\in \R^d$, $w_1,w_2\in \R^{d+1}$ that
\begin{align}\label{b04bc}
\max\!\left\{
\lvert g^d(x)\rvert, \lvert {g}_\varepsilon^d(x)\rvert\right\}\leq V_d(T,x),\quad
\max\!\left\{
\lvert Tf^d(0)\rvert, \lvert T{f}^d_\varepsilon(0)\rvert\right\}\leq V_d(t,x),
\end{align}
\begin{align}
\max\!\left\{
\lvert
f^d(w_1)-f^d(w_2)\rvert,
\lvert
{f}^d_\varepsilon(w_1)-{f}^d_\varepsilon(w_2)\rvert\right\}\leq \sum_{\nu=0}^{d}\left[
L_\nu^d\Lambda^d_\nu(T)
\lvert\pr^d_\nu(w_1-w_2) \rvert\right]
,\label{c01c}
\end{align}
\begin{align}
\max \!\left\{
\lvert g^d(x)-g^d(y)\rvert,
\lvert {g}^d_\varepsilon(x)-{g}^d_\varepsilon(y)\rvert
\right\}
\leq \frac{V_d(T,x)+V_d(T,y)}{2}
\frac{\lVert x-y\rVert}{\sqrt{T}},\label{b19c}
\end{align}
\begin{align}\label{k04c}
\left\lvert
g^d(x)-{g}^d_\varepsilon(x)\right\rvert\leq \varepsilon cd^cV_d(T,x).
\end{align}
\cref{b37} (applied for every $d\in \N$, $\varepsilon\in (0,1)$
with $\exponentV \gets 8$, $\exponentZ\gets 8$, $\exponentX\gets 8$, 
$(L_i)_{i\in [0,d]\cap \Z}\gets(L_i^d)_{i\in [0,d]\cap \Z}$,
$(\Lambda_i)_{i\in [0,d]\cap \Z}\gets(\Lambda_i^d)_{i\in [0,d]\cap \Z}$, $\pr\gets \pr^d$,
$f\gets f^d$,
$\tilde{f}\gets f^d_\varepsilon$,
$g\gets g^d$,
$\tilde{g}\gets g^d_\varepsilon$, $V\gets V_d$, $X\gets X^{d,0}$, $Z\gets Z^{d,0}$
in the notation of \cref{b37}), \eqref{b04bc}, \eqref{c01c}, \eqref{d05}, \eqref{d06}, \eqref{b19c}, \eqref{d07}, \eqref{d01}, \eqref{c74d}, \eqref{k03b}, and \eqref{k04c} show that the following items hold.

\begin{enumerate}[(a)]
\item\label{d09} For all 
$d\in \N$,
$\varepsilon\in (0,1)$ there exist  unique continuous functions $u^d, {u}^d_\varepsilon\colon [0,T)\times \R^d\to \R^{d+1}$ such that for all $t\in [0,T)$, $x\in \R^d$
we have that
\begin{align}
\max_{\nu\in [0,d]\cap\Z}\sup_{\tau\in [0,T), \xi\in \R^d}
\left[\Lambda^d_\nu(T-\tau)\frac{\lvert\pr^d_\nu(u^d(\tau,\xi))\rvert}{V_d(\tau ,\xi)}\right]<\infty,
\end{align}
\begin{align}
\max_{\nu\in [0,d]\cap\Z}\sup_{\tau\in [0,T), \xi\in \R^d}
\left[\Lambda^d_\nu(T-\tau)\frac{\lvert\pr^d_\nu({u}^d_\varepsilon(\tau,\xi))\rvert}{V_d(\tau ,\xi)}\right]<\infty,\label{d10}
\end{align}
\begin{align}\max_{\nu\in [0,d]\cap\Z}\left[
\E\!\left [\left\lvert g^d(X^{d,0,t,x}_{T} )\pr^d_\nu(Z^{d,0,t,x}_{T})\right\rvert \right] + \int_{t}^{T}
\E \!\left[\left\lvert
f^d(u^d(r,X^{d,0,t,x}_{r}))\pr^d_\nu(Z^{d,0,t,x}_{r})\right\rvert\right]dr\right]<\infty,
\end{align}
\begin{align}
\max_{\nu\in [0,d]\cap\Z}\left[
\E\!\left [\left\lvert {g}^d_\varepsilon(X^{d,0,t,x}_{T} )\pr_\nu^d(Z^{d,0,t,x}_{T})\right\rvert \right] + \int_{t}^{T}
\E \!\left[\left\lvert
{f}^d_\varepsilon({u}^d_\varepsilon(r,X^{d,0,t,x}_{r}))\pr_\nu(Z^{d,0,t,x}_{r})\right\rvert\right]dr\right]<\infty,
\end{align}
\begin{align}
u^d(t,x)=\E\!\left [g^d(X^{d,0,t,x}_{T} )Z^{d,0,t,x}_{T} \right] + \int_{t}^{T}
\E \!\left[
f^d(u^d(r,X^{d,0,t,x}_{r}))Z^{d,0,t,x}_{r}\right]dr,\label{d11}
\end{align}
and
\begin{align}
{u}^d_\varepsilon (t,x)=\E\!\left [{g}^d_\varepsilon(X^{d,0,t,x}_{T} )Z^{d,0,t,x}_{T} \right] + \int_{t}^{T}
\E \!\left[
{f}_\varepsilon^d({u}^d_\varepsilon(r,X^{d,0,t,x}_{r}))Z^{d,0,t,x}_{r}\right]dr.\label{d11b}
\end{align}

\item For all
 $t\in [0,T)$, $x\in \R^d$ we have that
\begin{align}
\max_{i\in [0,d]\cap\Z}
\left\lVert
\Lambda_{i}^d(T-t)
\pr_i^d({u}^d_\varepsilon(t,x) 
-u^d(t,x))\right\rVert_2
\leq 
10 c^2d^{2c}\varepsilon V_d^{3q+1}(t,x)\mathrm{B}(1-\tfrac{q}{2},\tfrac{1}{2}).\label{d13}
\end{align}
\end{enumerate}
Item \eqref{d09}, \cite[Theorem~6.9]{NW2023} (applied for every $d\in \N$ with $\mu^d\gets 0$, $\sigma^d \gets \mathrm{Id}_{\R^{d\times d}}$),  \eqref{c01bb}, and \eqref{b19b} 
imply that for all $d\in \N$ 
the function $u^d$ 
satisfies that
$v^d:=\pr_0^d(u^d)$ is the unique viscosity solution to the following 
semilinear PDE of parabolic type:
\begin{align}
&
\frac{\partial v^d}{\partial t}(t,x)
+\frac{1}{2}(\Delta v^d)(t,x)+f^d(t,x,v^d(t,x), (\nabla_xv^d)(t,x))
=0
\quad 
\forall\, t\in (0,T), x\in \R^d,
\\
&
v^d(T,x)=g^d(x) \quad \forall\, x\in \R^d,
\end{align}
$\nabla_xv^d=(\pr^d_1(u^d), \pr^d_2(u^d),\ldots, \pr^d_d(u^d))$,
and
\begin{align}
\max_{\nu\in [0,d]\cap\Z}\sup_{\tau\in [0,T), \xi\in \R^d}
\left[\Lambda^d_\nu(T-\tau)\frac{\lvert\pr^d_\nu(u^d(\tau,\xi))\rvert}{(1+\lVert \xi\rVert^2)^\frac{1}{2}}\right]<\infty.
\end{align}
This establishes \eqref{k36}.

Next,
\cite[Lemma 4.3]{NNW2023}
(applied for every $d\in \N$, $\varepsilon\in (0,1)$
with $\exponentV \gets 8$, $\exponentZ\gets 8$, $\exponentX\gets 8$, 
$(L_i)_{i\in [0,d]\cap \Z}\gets(L_i^d)_{i\in [0,d]\cap \Z}$,
$(\Lambda_i)_{i\in [0,d]\cap \Z}\gets(\Lambda_i^d)_{i\in [0,d]\cap \Z}$, $\pr\gets \pr^d$,
$f\gets f^d_\varepsilon$,
$g\gets g^d_\varepsilon$,
 $V\gets V_d$, $X\gets X^{d,0}$, $Z\gets Z^{d,0}$, $(U^\theta_{n,m})\gets (U^{d,\theta}_{n,m, \varepsilon})$, $q_1\gets 3$
in the notation of \cite[Lemma 4.3]{NNW2023})
demonstrate for all
$d,m,n\in \N$, $t\in [0,T)$, $x\in \R^d$, $\varepsilon \in (0,1)$ that
\begin{align}
\max_{\nu\in [0,d]\cap\Z}
\left\lVert\Lambda^d_\nu(T-t)\pr^d_\nu(
U^{d,0}_{n,m,\varepsilon}(t,x)-u^d_\varepsilon(t,x))
\right\rVert_2\leq e^\frac{m^3}{6}m^{-\frac{n}{2}}8^ne^{nc^2T}V_d^3(t,x).\label{d12}
\end{align}
Moreover, \eqref{d02} and \eqref{d03} show that there exists $\kappa\in (0,\infty)$ such that 
\begin{align}\label{d16}
\sqrt{d+1}
\left(
\int_{[0,1]^d}V^6_d(0,x)\, dx\right)^\frac{1}{2}+
\sqrt{d+1}
\left(
\int_{[0,1]^d}
V_d^{6q+2}(0,x)\, dx\right)^\frac{1}{2}\leq \kappa d^\kappa.
\end{align}
This, the triangle inequality, \eqref{d12}, and  \eqref{d13}
prove for all $d,n,m\in \N$, $\varepsilon\in (0,1)$
that
\begin{align}
&\left(
\int_{[0,1]^d}
\sum_{\nu=0}^{d}
\left\lVert\Lambda^d_\nu(T)\pr^d_\nu(
U^{d,0}_{n,m,\varepsilon}(0,x)-u^d(0,x))
\right\rVert_2^2dx\right)^\frac{1}{2}\nonumber\\
&\leq 
\left(
\int_{[0,1]^d}(d+1)
\max_{\nu\in [0,d]\cap\Z}
\left\lVert\Lambda^d_\nu(T)\pr^d_\nu(
U^{d,0}_{n,m,\varepsilon}(0,x)-u^d(0,x))
\right\rVert_2^2dx\right)^\frac{1}{2}\nonumber\\
&\leq \sqrt{d+1}
\left(
\int_{[0,1]^d}
\max_{\nu\in [0,d]\cap\Z}
\left\lVert\Lambda^d_\nu(T)\pr^d_\nu(
U^{d,0}_{n,m,\varepsilon}(0,x)-u^d_\varepsilon(0,x))
\right\rVert_2^2dx\right)^\frac{1}{2}\nonumber\\
&\quad +\sqrt{d+1}
\left(
\int_{[0,1]^d}
\max_{\nu\in [0,d]\cap\Z}
\left\lVert\Lambda^d_\nu(T)\pr^d_\nu(
u^d_\varepsilon(0,x)-u^d(0,x))
\right\rVert_2^2dx\right)^\frac{1}{2}\nonumber\\
&\leq  e^\frac{m^3}{6}m^{-\frac{n}{2}}8^ne^{nc^2T}
\sqrt{d+1}
\left(
\int_{[0,1]^d}V^6_d(0,x)\, dx\right)^\frac{1}{2}\nonumber\\&\quad 
+10 c^2d^{2c}\varepsilon \mathrm{B}(1-\tfrac{q}{2},\tfrac{1}{2})
\sqrt{d+1}\left(
\int_{[0,1]^d}
V_d^{6q+2}(0,x)\, dx\right)^\frac{1}{2}
\nonumber\\
&\leq 
\left(
e^\frac{m^3}{6}m^{-\frac{n}{2}}8^ne^{nc^2T}
+10 c^2d^{2c}\varepsilon \mathrm{B}(1-\tfrac{q}{2},\tfrac{1}{2})
\right)\kappa d^\kappa.\label{d17}
\end{align}
For the next step for every $d\in \N$, $\epsilon\in (0,1)$ let  $N_{d,\epsilon}\in \N$, $\varepsilon_{d,\epsilon}\in (0,1)$ satisfies that
\begin{align}
N_{d,\epsilon}=\inf\!\left\{n\in \N\cap[2,\infty)\colon 
e^\frac{n^3}{6}n^{-\frac{n}{2}}8^n e^{nc^2T}\kappa d^\kappa\leq \frac{\epsilon}{2}
\right\}\label{d14}
\end{align}
and 
\begin{align}\label{d15}
\varepsilon_{d,\epsilon}= \frac{\epsilon}{10 c^2d^{2c} \mathrm{B}(1-\tfrac{q}{2},\tfrac{1}{2})
\kappa d^\kappa}.
\end{align}Then Fubini's theorem and \eqref{d17} imply for all $d\in \N$, $\epsilon\in (0,1)$ that
\begin{align}
&\left(
\E \!\left[
\int_{[0,1]^d}\sum_{\nu=0}^{d}
\left\lvert\Lambda^d_\nu(T)\pr^d_\nu(
U^{d,0}_{N_{d,\epsilon},N_{d,\epsilon},\varepsilon_{d,\epsilon}}(0,x)-u^d(0,x))
\right\rvert^2dx\right]\right)^\frac{1}{2}
\nonumber\\
&
=\left(
\int_{[0,1]^d}
\sum_{\nu=0}^{d}
\E\!\left[
\left\lvert\Lambda^d_\nu(T)\pr^d_\nu(
U^{d,0}_{N_{d,\epsilon},N_{d,\epsilon},\varepsilon_{d,\epsilon}}(0,x)-u^d(0,x))
\right\rvert^2\right]dx\right) ^\frac{1}{2}\nonumber\\
&\leq
\left(
e^\frac{N_{d,\epsilon}^3}{6}N_{d,\epsilon}^{-\frac{N_{d,\epsilon}}{2}}8^{N_{d,\epsilon}}e^{N_{d,\epsilon}c^2T}
+10 c^2d^{2c}\varepsilon_{d,\epsilon} \mathrm{B}(1-\tfrac{q}{2},\tfrac{1}{2})
\right)\kappa d^\kappa\nonumber\\
&
\leq \frac{\epsilon}{2}+\frac{\epsilon}{2}=\epsilon.
\end{align}
This shows that for all $d\in \N$, $\epsilon\in (0,1)$ that there exists $\omega_{d,\epsilon}\in \Omega$ such that
\begin{align}\label{d24}
\int_{[0,1]^d}\sum_{\nu=0}^{d}
\left\lvert\Lambda^d_\nu(T)\pr^d_\nu(
U^{d,0}_{N_{d,\epsilon},N_{d,\epsilon},\varepsilon_{d,\epsilon}}(0,x)-u^d(0,x))
\right\rvert^2dx\leq \epsilon^2.
\end{align}
Next, 
\cref{c02} (applied for every 
$d\in \N$, $\epsilon\in (0,1)$ with $f\gets f^d_{\varepsilon_{d,\epsilon}}$,
$g\gets g^d_{\varepsilon_{d,\epsilon}}$,
$(U^\theta_{n,m})\gets (U^{d,\theta}_{n,m,\varepsilon_{d,\epsilon}})$, $\omega\gets \omega_{d,\epsilon}$ in the notation of \cref{c02})  imply for all $d\in \N$, $\epsilon\in (0,1)$ that
there exists $\Psi_{d,\epsilon}\in \bfN$
such that
\begin{align}\label{d18}
\dim (\calD(\Psi_{d,\epsilon}))= N_{d,\epsilon} (\dim (\calD (\Phi_{f^d_{\varepsilon_{d,\epsilon}}}))-1 )
+\dim (\calD(\Phi_{g^d_{\varepsilon_{d,\epsilon}}})),
\end{align}
\begin{align}\label{d19}
\supnorm{\calD(\Psi_{d,\epsilon})}\leq \max\!\left\{d+1,
\supnorm{\calD(\Phi_{g^d_{\varepsilon_{d,\epsilon}}})},
\supnorm{\calD(\Phi_{f^d_{\varepsilon_{d,\epsilon}}})} \right\}(4N_{d,\epsilon})^{N_{d,\epsilon}},
\end{align}
\begin{align}\label{d20}
U^{d,0}_{N_{d,\epsilon},N_{d,\epsilon},\varepsilon_{d,\epsilon}}(0,x,\omega_{d,\epsilon})= (\calR(\Psi_{d,\epsilon}))(x).
\end{align}
Furthermore, \eqref{d04} and the fact that 
$\forall\,\Phi\in \bfN \colon \max \{\dim(\calD(\Phi)), \supnorm{\calD(\Phi)}\}\leq \calP(\Phi)$ imply for all $d\in \N$, $\varepsilon\in (0,1)$ that
\begin{align}
\max \!\left\{
\dim (\calD(\Phi_{g^d_\varepsilon})),
\dim (\calD(\Phi_{f^d_\varepsilon})),\supnorm{\calD(\Phi_{g^d_\varepsilon})},
\supnorm{\calD(\Phi_{f^d_\varepsilon})}\right\}
\leq cd^c\varepsilon^{-c}.
\end{align}
Hence, we have for all $d\in \N$, $\epsilon\in (0,1)$ that
\begin{align}
\max\!\left\{
\dim (\calD(\Phi_{g^d_{\varepsilon_{d,\epsilon}}})),
\dim (\calD(\Phi_{f^d_{\varepsilon_{d,\epsilon}}})),\supnorm{\calD(\Phi_{g^d_{\varepsilon_{d,\epsilon}}})},
\supnorm{\calD(\Phi_{f^d_{\varepsilon_{d,\epsilon}}})}
\right\}
\leq cd^c\varepsilon^{-c}_{d,\epsilon}.
\end{align}
This, \eqref{d18},  \eqref{d19}, and the fact that $c\geq 2$ prove 
for all $d\in \N$, $\epsilon\in (0,1)$ that
\begin{align}
\dim (\calD(\Psi_{d,\epsilon}))\leq 2N_{d,\epsilon}cd^c\varepsilon_{d,\epsilon}^{-c}
\quad \text{and}\quad
\supnorm{\calD(\Psi_{d,\epsilon})}\leq cd^c\varepsilon_{d,\epsilon}^{-c} (4N_{d,\epsilon})^{N_{d,\epsilon}}.\label{d21}
\end{align}
Next, \eqref{d14} implies for all $d\in \N$, $\epsilon\in (0,1)$ that
\begin{align}
\epsilon&\leq 
e^\frac{(N_{d,\epsilon}-1)^3}{6}(N_{d,\epsilon}-1)^{-\frac{N_{d,\epsilon}-1}{2}}8^{N_{d,\epsilon}-1} e^{(N_{d,\epsilon}-1)c^2T}\kappa d^\kappa\nonumber\\
&\leq 
e^\frac{N_{d,\epsilon}^3}{6}(N_{d,\epsilon}-1)^{-\frac{N_{d,\epsilon}-1}{2}}8^{N_{d,\epsilon}} e^{N_{d,\epsilon}c^2T}\kappa d^\kappa.
\end{align}
This, the fact that
$\forall\,\Phi\in \bfN\colon \calP(\Phi)\leq 2 \dim (\calD(\Phi))
\supnorm{\calD(\Phi)}^2$,
\eqref{d21}, and \eqref{d15} show for all 
$d\in \N$, $\epsilon,\gamma\in (0,1)$ that
\begin{align}
\calP(\Psi_{d,\epsilon})&\leq 2 \dim (\calD(\Psi_{d,\epsilon}))
\supnorm{\calD(\Psi_{d,\epsilon})}^2\nonumber\\
&\leq 2N_{d,\epsilon}cd^c\varepsilon_{d,\epsilon}^{-c}
\left( cd^c\varepsilon_{d,\epsilon}^{-c} (4N_{d,\epsilon})^{N_{d,\epsilon}}\right)^2\nonumber\\
&=32 c^3 d^{3c}
\varepsilon_{d,\epsilon}^{-3c}N_{d,\epsilon}^{2N_{d,\epsilon}+1}
\nonumber\\
&\leq 32 c^3 d^{3c} \left(\frac{\epsilon}{10 c^2d^{2c} \mathrm{B}(1-\tfrac{q}{2},\tfrac{1}{2})
\kappa d^\kappa}\right)^{-3c}N_{d,\epsilon}^{2N_{d,\epsilon}+1}\nonumber\\
&=
32 c^3 d^{3c}\left(
10 c^2d^{2c} \mathrm{B}(1-\tfrac{q}{2},\tfrac{1}{2})
\kappa d^\kappa\right)^{3c}\epsilon^{-3c}
N_{d,\epsilon}^{2N_{d,\epsilon}+1}\epsilon^{4+\gamma}\epsilon^{-4-\gamma}\nonumber\\
&\leq 
32 c^3 d^{3c}\left(
10 c^2d^{2c} \mathrm{B}(1-\tfrac{q}{2},\tfrac{1}{2})
\kappa d^\kappa\right)^{3c}\epsilon^{-3c}
N_{d,\epsilon}^{2N_{d,\epsilon}+1}\nonumber\\
&\quad \left[
e^\frac{N_{d,\epsilon}^3}{6}(N_{d,\epsilon}-1)^{-\frac{N_{d,\epsilon}-1}{2}}8^{N_{d,\epsilon}} e^{N_{d,\epsilon}c^2T}\kappa d^\kappa\right]^{4+\gamma}\epsilon^{-4-\gamma}\nonumber\\
&=
32 c^3 d^{3c}\left(
10 c^2d^{2c} \mathrm{B}(1-\tfrac{q}{2},\tfrac{1}{2})
\kappa d^\kappa\right)^{3c}
(\kappa d^\kappa)^{4+\gamma}
\epsilon^{-3c-4-\gamma}
\sup_{n\in [2,\infty)\cap\Z}
\frac{n^{2n+1}(e^{\frac{n^3}{6}}8^ne^{nc^2T})^{4+\gamma}}{(n-1)^{\frac{n-1}{2}(\gamma+4) }}.\label{d22}
\end{align}
Next, note that for all $\gamma\in (0,1)$ we have that
\begin{align}
&
\sup_{n\in [2,\infty)\cap\Z}
\frac{n^{2n+1}(e^{\frac{n^3}{6}}8^ne^{nc^2T})^{4+\gamma}}{(n-1)^{\frac{n-1}{2}(\gamma+4) }}\nonumber\\
&\leq 
\sup_{n\in [2,\infty)\cap\Z}
\frac{n^{2n+1}(e^{\frac{n^3}{6}}8^ne^{nc^2T})^{4+\gamma}n^{\frac{\gamma+4}{2}}}{(n-1)^{\frac{n}{2}(\gamma+4) }}\nonumber\\
&=
\sup_{n\in [2,\infty)\cap\Z}\left[
\frac{n^{2n+1}(e^{\frac{n^3}{6}}8^ne^{nc^2T})^{4+\gamma}n^{\frac{\gamma+4}{2}}}{n^{\frac{n}{2}(\gamma+4) }}
\left(\frac{n}{n-1}\right)^{\frac{n}{2}(\gamma+4) }\right]\nonumber\\
&\leq 
\sup_{n\in [2,\infty)\cap\Z}\left[
\frac{n(e^{\frac{n^3}{6}}8^ne^{nc^2T})^{4+\gamma}n^{\frac{\gamma+4}{2}}}{n^{\frac{n\gamma}{2}}}
2^{\frac{n}{2}(\gamma+4) }\right]\nonumber\\
&
<\infty.
\end{align}
This and \eqref{d22} prove that there exists 
 $\eta\in (0,\infty)$ such that
for all $d\in \N$, $\epsilon\in (0,1)$ we have that
$
\calP(\Psi_{d,\epsilon})\leq  \eta d^\eta \epsilon^{-\eta}.
$
This, \eqref{d24}, and \eqref{d20} complete the proof of \cref{d23}.
\end{proof}

\section{Example}\label{s02}
In this section we give an example of
$(f^d_\varepsilon)_{d\in\N,\varepsilon\in (0,1)}$,
$(f^d)_{d\in\N}$,
$(g^d_\varepsilon)_{d\in\N,\varepsilon\in (0,1)}$,
$(g^d)_{d\in\N}$,
$(\Phi_{f^d_\varepsilon})_{d\in\N,\varepsilon\in (0,1)},(\Phi_{g^d_\varepsilon})_{d\in\N,\varepsilon\in (0,1)}$
that satisfy
the assumption of \cref{d23}, see \cref{e01} below. First of all, we need the following result.
\begin{lemma}\label{c05}
Assume \cref{m07b}. Let $d\in \N$. Let $f\colon \R^d\to\R$ satisfy for all $x=(x_i)_{i\in [1,d]\cap\Z}\in \R^d$ that $f(x)=\frac{1}{d}\sum_{i=1}^{d}x_i$. Then $f\in \calR (\{\Phi\in \bfN\colon \calD(\Phi)=(d,2d,1) \})$.
\end{lemma}
\begin{proof}
[Proof of \cref{c05}]We adapt the proof of \cite[Lemma~3.6]{CHW2022}.
Let $ w_1=\binom{1}{-1} $ and  let $ W_1 \in \R^{2d\times d} $, $ B_1\in \R^{2d} $, $W_2\in \R^{1\times 2d}$, $B_1\in \R $ satisfy
\begin{align} \begin{split} 
&
W_1=
\begin{pmatrix}
w_1 & 0 &\ldots&0 \\
0 & w_1&\ldots&0\\
\vdots&\vdots &\ddots&\vdots\\
0 &0 &\ldots&w_1
\end{pmatrix}
\in \R^{2d\times d},\quad B_1=
0_{\R^{2d}}
\in\R^{2d},
\\
&
W_2=\frac{1}{d}\left(w_1^\top,w_1^\top,\ldots,w_1^\top\right), \quad B_2=0.
\end{split}
\end{align}
Let $ \phi\in \mathbf{N} $,
  $ x^{0} =(x^{0}_1,...,x^{0}_d)\in \R^d $, $ x^{1},x^{2} \in \R^{2d} $ satisfy that 
$ \phi=((W_1,B_1),(W_2,B_2) )$,	$x^{1}=\mathbf{A}_{2d}(W_1x^{0}+B_1)
$ and $x^2= W_2x^1+B_2$.
Then $\calD(\phi) =(d,2d,1) $ and
\begin{align}
x^1=\mathbf{A}_{2d}\left  ((x_1^0,-x_1^0,...,x_d^0,-x_d^0)^\top\right )=\left ((x_1^0)^+,(-x_1^0)^+,...,(x_d^0)^+,(-x_d^0)^+\right )^\top,
\end{align}
where we write
$(\cdot )^+= \max \{\cdot ,0\}$. Hence, we have that
$x^2=W_2x^1+B_2=\frac{1}{d}\sum_{i=1}^{d}((x^0_i)^+-(-x^0_i)^{+})=\frac{1}{d}\sum_{i=1}^{d}x^0_i$. Thus, by definition we have
for all $x=(x_i)_{i\in [1,d]\cap\Z}\in \R^d$
 that $(\calR(\phi))(x)= \frac{1}{d}\sum_{i=1}^{d}x_i=f(x)$. This and the fact that $\calD(\phi) =(d,2d,1) $ complete the proof of \cref{c05}.
\end{proof}

\begin{example}\label{e01}Assume \cref{m07}.
For every $d\in \N$, $\varepsilon\in (0,1)$
let $f^d\in C(\R^{d+1},\R)$, $g^d\in C(\R^d,\R)$ satisfy for all
$w=(w_i)_{i\in [0,d]\cap\Z}$,
$x=(x_i)_{i\in [1,d]\cap\Z}$ that
$f^d(w)=\sin (\frac{1}{d+1}\sum_{i=0}^d w_i)$ and
$g^d(w)=\sin (\frac{1}{d}\sum_{i=1}^d x_i)$. 
Then there exist
$(f^d_\varepsilon)_{d\in\N,\varepsilon\in (0,1)}\subseteq C(\R^{d+1},\R)$,
$(g^d_\varepsilon)_{d\in\N,\varepsilon\in (0,1)}\subseteq C(\R^d,\R)$,
$(\Phi_{f^d_\varepsilon})_{d\in\N,\varepsilon\in (0,1)},(\Phi_{g^d_\varepsilon})_{d\in\N,\varepsilon\in (0,1)}\subseteq \bfN$ which satisfy that
$\calR(\Phi_{f^d_\varepsilon})=f^d_\varepsilon$ and
$\calR(\Phi_{g^d_\varepsilon})=g^d_\varepsilon$
and that \eqref{b04bb}--\eqref{d04} hold for 
$T=1$, $\beta=2$, $q=1.5$ and for
some $c\in [2,\infty)$.
\end{example}
\begin{proof}[Proof of \cref{e01}] Throughout this proof 
let \cref{m07b} be given.
First, \cite[Corollary 3.13]{HJKN2020a} (applied for every
$\varepsilon\in (0,1)$ with
$\epsilon\gets \varepsilon$,
$L\gets 1$,
$q\gets 1.5$,
$f\gets \sin (\cdot) $ in the notation of 
 \cite[Corollary 3.13]{HJKN2020a}) and the fact that
$\forall\,x,y\in \R\colon \lvert\sin(x) -\sin(y)\rvert\leq \lvert x-y\rvert$ ensure  that
there exist $\kappa\in [1,\infty)$,
$(\gamma_\varepsilon)_{\varepsilon\in (0,1)}\subseteq C(\R,\R)$, $(\Phi_{\gamma_\varepsilon})_{\varepsilon\in (0,1)}\subseteq \bfN$ such that for all $\varepsilon\in (0,1)$, $x,y\in \R$ we have that
\begin{align}
\left
\lvert
\gamma_\varepsilon(x)-\gamma_\varepsilon(y)\right\rvert
\leq \lvert x-y\rvert,
\quad 
\left\lvert\sin(x)-\gamma_\varepsilon(x)\right\rvert
\leq \varepsilon\left(1+\lvert x\rvert^{1.5}\right),\label{x01}
\end{align}
$\calD(\Phi_{\gamma_\varepsilon})\in \N^3$, 
$\supnorm{\calD(\Phi_{\gamma_\varepsilon})}\leq \kappa\varepsilon^{-3}$, and 
$\calR(\Phi_{\gamma_\varepsilon})=\gamma_\varepsilon$.
For every $d\in \N$, $\varepsilon\in (0,1)$
let $f^d_\varepsilon\in C(\R^{d+1},\R)$, $g^d_\varepsilon\in C(\R^d,\R)$ satisfy for all
$w=(w_i)_{i\in [0,d]\cap\Z}$,
$x=(x_i)_{i\in [1,d]\cap\Z}$ that
$f^d_\varepsilon(w)=\gamma_\varepsilon (\frac{1}{d+1}\sum_{i=0}^d w_i)$
 and
$g^d_\varepsilon(w)=\gamma_\varepsilon (\frac{1}{d}\sum_{i=1}^d x_i)$.
Next, 
\cref{c05} shows  for all $d\in \N$ that
$
(\R^d\ni x=(x_i)_{i\in [1,d]\cap\Z}\mapsto \frac{1}{d}\sum_{i=1}^{d}x_i\in\R)\in \calR (\{\Phi\in \bfN\colon \calD(\Phi)= (d,2d,1) \})
$. This,  \cref{m11b}, and the definition of 
$(f^d_\varepsilon)_{d\in\N,\varepsilon\in (0,1)}$,
$(g^d_\varepsilon)_{d\in\N,\varepsilon\in (0,1)}$, $(\gamma_\varepsilon)_{\varepsilon\in (0,1)}$,
$(\Phi_{\gamma_\varepsilon})_{\varepsilon\in (0,1)}$
show that
there exist $(\Phi_{f^d_\varepsilon})_{d\in\N,\varepsilon\in (0,1)},(\Phi_{g^d_\varepsilon})_{d\in\N,\varepsilon\in (0,1)}\subseteq \bfN$
such that for all 
$d\in\N$, $\varepsilon\in (0,1)$ we have that
$\calR(\Phi_{f^d_\varepsilon})=f^d_\varepsilon$,
$\calR(\Phi_{g^d_\varepsilon})=g^d_\varepsilon$,
$\calD(\Phi_{f^d_\varepsilon}) = \calD(\Phi_{\gamma_\varepsilon})\odot(d+1,2d+2,1) $,
$\calD(\Phi_{g^d_\varepsilon}) =
\calD(\Phi_{\gamma_\varepsilon})\odot(d,2d,1)
$.
This, the definition of $\odot$, and the fact that
$\forall\,\varepsilon\in (0,1)\colon\supnorm{\calD(\Phi_{\gamma_\varepsilon})}\leq \kappa\varepsilon^{-3}$
 show that
$\dim (\calD(\Phi_{f^d_\varepsilon}) )=
\dim (\calD(\Phi_{g^d_\varepsilon}) )=5
$
and 
$
\max \!\left\{
\supnorm{\calD(\Phi_{f^d_\varepsilon}) },
\supnorm{\calD(\Phi_{g^d_\varepsilon}) }
\right\}
\leq \supnorm{\calD(\Phi_{\gamma_\varepsilon})}+4d
\leq \kappa\varepsilon^{-3}+4d\leq 8\kappa d\varepsilon^{-3}
$.
This and the fact that
$\forall\,\Phi\in \bfN\colon \calP(\Phi)\leq 2 \dim (\calD(\Phi))
\supnorm{\calD(\Phi)}^2$ show that
\begin{align}
\max \!\left\{
{\calP(\Phi_{f^d_\varepsilon}) },
{\calP(\Phi_{g^d_\varepsilon}) }
\right\}
\leq 2\cdot 5(8\kappa d\varepsilon^{-3})^2
=640\kappa^2d^2\varepsilon^{-6}
.
\end{align} Furthermore, the definition of 
$(f^d_\varepsilon)_{d\in\N,\varepsilon\in (0,1)}$, the Lipschitz condition of $\gamma_\varepsilon$, $\varepsilon\in (0,1)$, and the triangle inequality show for all 
$d\in \N$, $\varepsilon\in (0,1)$, $w^1=(w^1_i)_{i\in [0,d]\cap\Z},
w^2=(w^2_i)_{i\in [0,d]\cap\Z}
\in \R^{d+1}
$
that
\begin{align} \begin{split} 
\left\lvert
f^d_\varepsilon(w^1)-
f^d_\varepsilon(w^2)\right\rvert
&
=
\left\lvert
\gamma_\varepsilon\left(\frac{1}{d+1}\sum_{i=0}^d w^1_i\right)
-
\gamma_\varepsilon\left(\frac{1}{d+1}\sum_{i=0}^d w^2_i\right)
\right\rvert\\
&\leq \left\lvert\frac{1}{d+1}\sum_{i=0}^d (w^1_i-w^2_i)\right\rvert\\
&\leq \frac{1}{d+1}\sum_{i=0}^d\left\lvert w^1_i-w^2_i\right\rvert.
\end{split}
\end{align}
A similar argument and the Cauchy-Schwarz inequality show for all 
$d\in \N$, $\varepsilon\in (0,1)$, $x=(x_i)_{i\in [1,d]\cap\Z},
y=(y_i)_{i\in [1,d]\cap\Z}
\in \R^{d}
$
that
\begin{align}
\left\lvert
g^d_\varepsilon(x)-
g^d_\varepsilon(y)\right\rvert
\leq  \frac{1}{d}\sum_{i=1}^d\left\lvert x_i-y_i\right\rvert
\leq \left(
 \frac{1}{d}\sum_{i=1}^d\left\lvert x_i-y_i\right\rvert^2\right)^{\frac{1}{2}}= \frac{1}{\sqrt{d}}\lVert x-y\rVert.
\end{align}
Next,  the definition of 
$(f^d_\varepsilon)_{d\in\N,\varepsilon\in (0,1)}$,
$(f^d)_{d\in\N}$, \eqref{x01}, and Jensen's inequality show for all
$d\in \N$, $\varepsilon\in (0,1)$, $w=(w_i)_{i\in [0,d]\cap\Z}\in \R^{d+1}$
that
\begin{align} \begin{split} 
\left\lvert
f^d_\varepsilon(w)-f^d(w)\right\rvert
&=
\left\lvert
\gamma_\varepsilon\left(\frac{1}{d+1}\sum_{i=0}^d w_i\right)-\sin \left(\frac{1}{d+1}\sum_{i=0}^d w_i\right)\right\rvert\\
&
\leq \varepsilon\left(1+\left\lvert
\frac{1}{d+1}\sum_{i=0}^d w_i
\right\rvert^{1.5}\right)\\
&\leq 
\varepsilon\left(1+
\frac{1}{d+1}\sum_{i=0}^d\left\lvert w_i\right\rvert^{1.5}
\right).\end{split}
\end{align}
A similar argument and Jensen's inequality
show for all
$d\in \N$, $\varepsilon\in (0,1)$, $x=(x_i)_{i\in [1,d]\cap\Z}\in \R^{d}$
that
\begin{align}
\left\lvert
g^d_\varepsilon(x)-g^d(x)\right\rvert\leq 
\varepsilon\left(1+\left\lvert
\frac{1}{d}\sum_{i=1}^d x_i
\right\rvert^{1.5}\right)
\leq 
\varepsilon\left(1+\left\lvert
\frac{1}{d}\sum_{i=1}^d \lvert x_i\rvert^2
\right\rvert^{\frac{3}{4}}\right)\leq \varepsilon(1+\lVert x\rVert^2).
\end{align}
It is now easy to see that
$(f^d_\varepsilon)_{d\in\N,\varepsilon\in (0,1)}$,
$(f^d)_{d\in\N}$,
$(g^d_\varepsilon)_{d\in\N,\varepsilon\in (0,1)}$,
$(g^d)_{d\in\N}$,
$(\Phi_{f^d_\varepsilon})_{d\in\N,\varepsilon\in (0,1)},(\Phi_{g^d_\varepsilon})_{d\in\N,\varepsilon\in (0,1)}$ satisfy   that
$\calR(\Phi_{f^d_\varepsilon})=f^d_\varepsilon$ and
$\calR(\Phi_{g^d_\varepsilon})=g^d_\varepsilon$
and that \eqref{b04bb}--\eqref{d04} hold for 
$T=1$, $\beta=2$, $q=1.5$ and for
some $c\in [2,\infty)$. 
\end{proof}
\section{Conclusion}
We have rigorously proven that DNNs with ReLU activation function can overcome the curse dimensionality in the approximation of heat equations \eqref{c36} with \emph{gradient-dependent} nonlinearities in the sense that the required number of parameters in the deep neural networks increases at most polynomially in both the dimension $d$ of \eqref{c36} and the reciprocal of the prescribed accuracy $\epsilon$.

One natural question for future research would be to analyze whether one can extend \cref{d23} for general non-polynomial activation function. As first step, one could extend \cref{d23} to leaky ReLU and softplus activation as done in \cite{AJK+2023} for nonlinear Kolmogorov equations without gradient-dependent nonlinearities. To this end one would need to show, i.a., that the identity in $\R$ can also be represented by DNNs with arbitrary lengths (cf.\ \cref{b03} and, e.g., \eqref{h09b}), which is doable by combining \cite[Section 3.2]{AJK+2023} and \cite[Lemmas 2.2.9--10]{JKvW2023}. However, for other activation functions, it is not clear at all if \cref{d23} can be extended accordingly, and we leave this for future research.

Moreover, in future research one could also try to extend \cref{d23} by considering general $L^p$-norms, $p\in [2,\infty)$, instead of only $p=2$, as analyzed in \cite{AJK+2023}  in the setting described above.
Furthermore, one may want to try to extend \eqref{c36} from semilinear heat equations to semilinear PDEs with general drift and diffusion coefficients,  as analyzed in \cite{CHW2022} in a setting  with corresponding nonlinear PDEs without gradient-dependent nonlinearities.

{
\bibliographystyle{acm}
\bibliography{References}
}

\end{document}